\documentclass[12pt, reqno]{amsart}

\usepackage{hyperref}

\usepackage{epsfig}
\usepackage[all]{xy}
\usepackage{amsmath,amscd,amsthm,amssymb,amsfonts}

\theoremstyle{plain}

\newtheorem{Prop}{Proposition}[section]
\newtheorem{Lem}[Prop]{Lemma}

\newtheorem{Thm}{Theorem}[section]

\theoremstyle{definition}
\newtheorem*{Def}{Definition}

\theoremstyle{remark}
\newtheorem{Ex}[Prop]{Example}
\newtheorem{Remark}[Prop]{Remark}

\newcommand{\renumi}{\renewcommand\labelenumi}

\newcommand{\condition}[1]{\textup{({#1})}}
\newcommand{\conditionenumi}{\condition{\roman{enumi}}}
\newcommand{\statement}[1]{\textup{({#1})}}
\newcommand{\statementenumi}{\statement{\arabic{enumi}}}

\newcommand{\arw}[1]{\ar@{-}[#1]\ar@{>}[#1]}
\newcommand{\arwex}[1]{\ar@<0.40ex>[#1]}
\newcommand{\arwexd}[1]{\ar@<-0.40ex>[#1]}
\newcommand{\arweq}[1]{\ar@{=>}[#1]}
\newcommand{\arwbr}[1]{\ar@{-->}[#1]}
\newcommand{\arwstrut}[1]{\ar@{}[#1]}
\newcommand{\aru}[2]{\arw{#1}^-{\displaystyle{#2}}}
\newcommand{\ard}[2]{\arw{#1}_-{\displaystyle{#2}}}

\newcommand{\arsub}{\ar@{}[r]|{\subset}}
\newcommand{\arcap}{\ar@{}[d]|{\bigcap}}
\newcommand{\arcup}{\ar@{}[u]|{\bigcup}}
\newcommand{\ariso}{\ar@{}[r]|{\cong}}

\newcommand{\arisod}{\ard{d}{\rotatebox{90}{$\sim$}}}

\newcommand{\squarerd}[8]{\xymatrix{{#1} \aru{r}{#2} \ard{d}{#3}& {#4} \aru{d}{#5}\\ {#6} \ard{r}{#7}& {#8}}}
\newcommand{\squarerdcomm}[8]{\xymatrix{{#1} \aru{r}{#2} \ard{d}{#3}& {#4} \aru{d}{#5}\\ {#6} \ard{r}{#7}& {#8} \ar@{}[ul]|{\circlearrowleft}}}

\newcommand{\upperx}[1]{{#1}_{\gqx}}
\newcommand{\lowerx}[1]{{#1}_{\lqx}}
\newcommand{\lowerbx}[1]{\operatorname{Low}_{#1}}

\newcommand{\px}{\oplus}
\newcommand{\dx}{\odot}
\newcommand{\sx}{\oslash}
\newcommand{\mx}{- \infty}

\newcommand{\infx}[3]{\inf_{#1} \{{#2},{#3} \}}

\newcommand{\bigpx}{\bigoplus}

\newcommand{\bigdirdx}{\prod}
\newcommand{\che}[1]{{#1}^{\vee}}

\newcommand{\chepx}{\operatorname{\overset{\vee}{\oplus}}}
\newcommand{\chedx}{\operatorname{\overset{\vee}{\odot}}}
\newcommand{\cheloc}[1]{\operatorname{Loc}({#1})}

\newcommand{\lqx}{\leq}
\newcommand{\gqx}{\geq}
\newcommand{\nlqx}{\nleq}

\newcommand{\ltx}{<}
\newcommand{\gtx}{>}
\newcommand{\rrt}{\mathbb{T}}

\newcommand{\pairx}[2]{\langle{#1},{#2} \rangle}

\newcommand{\rii}{\mathbb{I}}
\newcommand{\nn}{\mathbb{N}}
\newcommand{\pp}{\mathbb{P}}

\newcommand{\rr}{\mathbb{R}}

\newcommand{\zz}{\mathbb{Z}}
\newcommand{\cc}{\mathbb{C}}

\newcommand{\e}{\varepsilon}

\newcommand{\ph}{\varphi}

\newcommand{\aaa}{\mathcal{A}}

\newcommand{\ooo}{\mathcal{O}}

\newcommand{\rmap}{\longrightarrow}
\newcommand{\then}{\Rightarrow}
\newcommand{\iso}{\cong}
\newcommand{\setmid}{\,|\,}
\newcommand{\sethyphen}{\cdot}
\newcommand{\inv}[1]{{#1}^{-1}}
\newcommand{\res}[1]{|_{#1}}
\newcommand{\setexists}{\mathstrut^{\exists}}
\newcommand{\setforall}{\mathstrut^{\forall}}

\newcommand{\abso}[1]{|{#1}|}
\newcommand{\absow}[1]{||{#1}||}

\newcommand{\xxsup}{\operatorname*{sup}}
\newcommand{\xxinf}{\operatorname*{inf}}

\newcommand{\xxDiv}{\operatorname{Div}}

\newcommand{\xxHom}{\operatorname{Hom}}

\newcommand{\xxid}{\operatorname{id}}

\newcommand{\xxdim}{\operatorname{dim}}

\newcommand{\xxdet}{\operatorname{det}}
\newcommand{\xxord}{\operatorname{ord}}

\newcommand{\xxval}{\operatorname{val}}
\newcommand{\xxexp}{\operatorname{exp}}

\begin {document}

\title {Generators of modules in tropical geometry}
\author {Shuhei Yoshitomi}
\date{}



\maketitle

\mark{}


\tableofcontents


\section {Introduction}

\subsection {Results}

A tropical curve is a geometric object over the tropical semifield of real numbers $\rrt= (\rr \cup \{ \mx \}, \px, \dx)$, where the addition $\px$ is the max-operation in the real field $\rr$, and the multiplication $\dx$ is the addition of $\rr$. For a tropical curve $C$ and a divisor $D$ on $C$, the set $M= H^0(C, \ooo_{ C}( D))$ of the sections of $D$ has the structure of a $\rrt$-module that is defined as follows.
\par
A $\rrt$-module $M$ is defined as a module over a semifield. $(M, \px, \dx, \mx)$ is said to be a $\rrt$-module if $(M, \px, \mx)$ is a tropical semigroup, and $\dx$ is an additive semigroup action on $M$ by $\rrt$. A tropical semigroup is a commutative semigroup with unity such that any element $v$ satisfies the idempotent condition $v \px v= v$.
\par
A $\rrt$-module $M$ is analogous to a module over a field. A subset $S \subset M$ is said to be a basis if it is a minimal system of generators. But the number of elements of a basis of $M$ is not necessarily equal to the topological dimension of it. We introduce straight $\rrt$-modules in section \ref{Sec:FB6}. This class is a generalization of lattice-preserving submodules of the free $\rrt$-module $\rrt^n$, where a lattice-preserving submodule is a submodule preserving the infimum of any two elements with respect to the canonical partial order relation on $\rrt^n$.

\begin {Thm}
Let $M$ be a finitely generated straight submodule of the free $\rrt$-module $\rrt^n$. Then $M$ is generated by $n$ elements.
\end {Thm}

We have four corollaries (Theorem \ref{Prop:FE1}, \ref{Prop:FE2}, \ref{Prop:FE3}, \ref{Prop:FE4}). The semifield $\rrt$ is generalized to a quasi-complete totally ordered rational tropical semifield $k$. We find a sufficient condition to the existence of a left-inversion of an injective homomorphism of $k$-modules (Theorem \ref{Prop:FE1}). The dimension of a straight reflexive $k$-module is defined to be the number of elements of a basis. We show the inequality $\xxdim( M) \leq \xxdim( N)$ for a pair of straight reflexive $k$-modules $M \subset N$ (Theorem \ref{Prop:FE2}). We show that a finitely generated straight pre-reflexive $k$-module is reflexive (Theorem \ref{Prop:FE3}). Also we consider finiteness of a submodule of a $k$-module (Theorem \ref{Prop:FE4}). The proofs are given in section \ref{Sec:FB1}.
\par
This result has an application to polytopes in a tropical projective space $\rrt \pp^n$. By Joswig and Kulas \cite{Jos0801}, a polytrope (it means a polytope in $\rrt \pp^n$ that is real convex) is a tropical simplex, and therefore it is the tropically convex hull of at most $n+1$ points. We show a generalization of this result (Theorem \ref{Prop:FE5}). A polytope $P$ is the tropically convex hull of at most $n+1$ points if the corresponding submodule $M \subset \rrt^{ n+1}$ is straight reflexive. Also $M$ is straight reflexive if $P$ is a polytrope.
\par
Also we have an application to tropical curves. A Riemann-Roch theorem for tropical curves is proved by Gathmann and Kerber \cite{Gath0612}. This theorem states an equality for an invariant $r( D)$ of the divisor. We see that $r( D)$ is not an invariant of the $\rrt$-module $M= H^0(C, \ooo_{ C}( D))$ (Example \ref{Ex:FB1}), and show the inequality $r( D) \leq \xxdim( M)- 1$ (Theorem \ref{Prop:FE7}).

\subsection {Background}

Surveys of tropical mathematics are found in \cite{Litvinov0501}, \cite{Tsai1011}. Early studies of tropical curves are found in \cite{Gath0612}, \cite{Mik}, \cite{Mik2}. Tropical varieties are introduced as follows. Let $K= \cc[[ \rr]]$ be the group algebra of power series defined by the group $\rr$. We have a multiplicative seminorm
\[
\absow{ \cdot} \colon K \rightarrow \rr_{ \gqx 0}
\]
defined by
\[
\absow{ x}= \xxexp( - \xxval( x)),
\]
where $\xxval$ means the canonical valuation on $K$. This seminorm induces the amoeba map
\[
\aaa \colon ( K^{ \times})^n \rmap \rr^n
\]
defined by
\[
\aaa( x_1, \ldots, x_n)=( \log \absow{ x_1}, \ldots, \log \absow{ x_n}).
\]
The image $\aaa( V)$ of a variety $V$ in the algebraic torus $( K^{ \times})^n$ is said to be a tropical variety in the tropical torus $\rr^n$.
\par
Tropical algebra is introduced by the map
\[
\pi \colon K \rmap \rr \cup \{ \mx \}
\]
defined by
\[
\pi( x)= \log \absow{ x}.
\]
This map induces a hyperfield homomorphism
\[
\pi \colon K \rmap X,
\]
where $X$ is the tropical hyperfield with underlying set $\rr \cup \{ \mx \}$, introduced in \cite{Viro1006}. The power set $2^{ X}$ is a semiring with operations induced by multi-operations of $X$.
\par
Now we have the lower-saturation map
\[
\nu \colon X \rmap 2^{ X}
\]
defined by
\[
\nu( a)= \{ c \in X \setmid c \lqx a \}.
\]
The power set $2^{ X}$ has a subsemiring
\[
\rii = X \cup \nu( X),
\]
which is isomorphic to Izhakian's extended tropical semiring introduced in \cite{Izhakian0511}. The lower-saturation map $\nu$ means the ghost map in \cite{Izhakian0511}. The image $\nu( X)$ means the ghost part, which is isomorphic to the tropical semifield of real numbers $(\rr \cup \{ \mx \}, \px, \dx)$, where operations are defined as follows.
\[
a \px b= \max \{a, b \},
\]
\[
a \dx b= a+ b.
\]
In this paper, the symbol $\rrt$ means the tropical semifield of real numbers. Under the identification $\rrt= \nu( X)$, the canonical homomorphism $\nu \colon \rii \rightarrow \rrt$ is the lower-saturation map.

\vspace{1em}

Section \ref{Sec:FB6} contains definitions and theorems. Section \ref{Sec:FB2} and \ref{Sec:FB3} contain foundation of tropical modules, and the proof of Theorem \ref{Prop:FE1}, \ref{Prop:FE2}, \ref{Prop:FE3}, \ref{Prop:FE4}, and \ref{Prop:FE5}. Section \ref{Sec:FB4} and \ref{Sec:FB5} contain foundation of tropical matrices and tropical curves, and the proof of Theorem \ref{Prop:FE7}. Section \ref{Sec:B1} is an appendix for tropical plane curves.

\vspace{1em}

Acknowledgements. The author thanks Professor Yujiro Kawamata for helpful advice.

\section {Definitions and theorems} \label {Sec:FB6}

A {\em semigroup} $(M, \px)$ is a set $M$ with an associative operation $\px$.

\begin {Def}
$(M, \px, \mx)$ is a {\em tropical semigroup} if it satisfies the following axioms.
\begin {enumerate}
\renumi{\conditionenumi}
\item $(M, \px)$ is a semigroup.
\item $v \px w= w \px v$.
\item $v \px \mx= v$.
\item $v \px v= v$.
\end {enumerate}
The element $\mx$ is called the zero element of $M$.
\end {Def}

There is a unique partial order relation `$\lqx$' on $M$ such that for any $v, w \in M$ it implies
\[
\xxsup \{v, w \}= v \px w.
\]
The proof is given in section \ref{Sec:FA4}.

\begin {Def}
A tropical semigroup $M$ is {\em quasi-complete} if any non-empty subset $S \subset M$ admits the infimum $\inf(S)$ (i.e. it admits the maximum element of the lower-bounds of $S$).
\end {Def}

\begin {Def}
$(A, \px, \dx, \mx, 0)$ is a {\em tropical semiring} if it satisfies the following axioms.
\begin {enumerate}
\renumi{\conditionenumi}
\item $(A, \px, \mx)$ is a tropical semigroup.
\item $(A, \dx)$ is a semigroup.
\item $a \dx b= b \dx a$.
\item $a \dx (b \px c)= a \dx b \px a \dx c$.
\item $a \dx 0= a$.
\item $a \dx \mx= \mx$.
\end {enumerate}
The element $\mx$ is called the zero element of $A$. The element $0$ is called the unity of $A$.
\end {Def}

\begin {Def}
$(k, \px, \dx, \mx, 0)$ is a {\em tropical semifield} if it satisfies the following axioms.
\begin {enumerate}
\renumi{\conditionenumi}
\item $(k, \px, \dx, \mx, 0)$ is a tropical semiring.
\item For any $a \in k \setminus \{ \mx \}$ there is an element $\sx a \in k$ such that $a \dx( \sx a)= 0$.
\end {enumerate}
\end {Def}

\begin {Def}
A tropical semifield $k$ is {\em rational} if it satisfies the following conditions.
\begin {enumerate}
\renumi{\conditionenumi}
\item $a \in k$, $m \in \nn$ $\then$ $\setexists b \in k$, $a= b^{ \dx m}$.
\item $k$ has no maximum element.
\end {enumerate}
\end {Def}

The {\em tropical semifield of real numbers} $(\rrt, \px, \dx, \mx, 0)$ is the set
\[
\rrt= \rr \cup \{ \mx \}
\]
equipped with addition
\[
a \px b= \max \{a, b \}
\]
and multiplication
\[
a \dx b= a+ b
\]
and zero element $\mx$ and unity $0$. $\rrt$ is a quasi-complete totally ordered rational tropical semifield.
\par
Let $k$ be a quasi-complete totally ordered rational tropical semifield.

\begin {Def}
$(M, \px, \dx, \mx)$ is a {\em $k$-module} if it satisfies the following axioms.
\begin {enumerate}
\renumi{\conditionenumi}
\item $(M, \px, \mx)$ is a tropical semigroup.
\item $\dx$ is a semigroup action $k \times M \ni (a, v) \mapsto a \dx v \in M$, i.e.
\begin {enumerate}
\item[i)] $(a \dx b) \dx v= a \dx (b \dx v)$.
\item[ii)] $0 \dx v= v$.
\end {enumerate}
\item $(a \px b) \dx v= (a \dx v) \px (b \dx v)$.
\item $a \dx (v \px w)= (a \dx v) \px (a \dx w)$.
\item $\mx \dx v= \mx$.
\item $a \dx \mx= \mx$.
\end {enumerate}
\end {Def}

\begin {Def}
A {\em homomorphism} $\alpha \colon M \rightarrow N$ {\em of $k$-modules} is a map with the following conditions.
\begin {enumerate}
\renumi{\conditionenumi}
\item $\alpha( \mx)= \mx$.
\item $\alpha( v \px w)= \alpha(v) \px \alpha(w)$.
\item $\alpha( a \dx v)= a \dx \alpha( v)$.
\end {enumerate}
\end {Def}

Let $\xxHom( M, N)$ denote the $k$-module of homomorphisms from $M$ to $N$.
\par
The dual module $\che{ M}$ is defined by $\che{ M}= \xxHom( M, k)$. We have the pairing map $\pairx{ \sethyphen}{ \sethyphen} \colon M \times \che{ M} \rightarrow k$ defined by
\[
\pairx{ v}{ \xi}= \xi( v).
\]

\begin {Def}
$M$ is {\em pre-reflexive} if the homomorphism $\iota_{ M} \colon M \rightarrow \che{( \che{ M})}$ is injective. $M$ is {\em reflexive} if $\iota_{ M}$ is an isomorphism.
\end {Def}

\begin {Def}
A $k$-module $M$ is {\em straight} if it is a finitely distributive ordered lattice, i.e. it satisfies the following conditions.
\begin {enumerate}
\renumi{\conditionenumi}
\item Any two elements $v, w \in M$ admit the infimum $\infx{ M}{ v}{ w}$.
\item $v_1, v_2, w \in M \then \infx{ M}{ v_1 \px v_2}{ w}= \infx{ M}{ v_1}{ w} \px \infx{ M}{ v_2}{ w}$.
\item $v_1, v_2, w \in M \then \infx{ M}{ v_1}{ v_2} \px w= \infx{ M}{ v_1 \px w}{ v_2 \px w}$.
\end {enumerate}
\end {Def}

\begin {Def}
A homomorphism $\alpha \colon M \rightarrow N$ is {\em lightly surjective} if for any $w \in N$ there is $v \in M$ such that $w \lqx \alpha( v)$.
\end {Def}

A homomorphism $\beta \colon N \rightarrow M$ is said to be a left-inversion of $\alpha$ if $\beta \circ \alpha= \xxid_{ M}$.

\begin {Thm} \label {Prop:FE1}
Let $\alpha \colon M \rightarrow N$ be an injective lightly surjective homomorphism of $k$-modules such that $M$ is straight reflexive. Then $\alpha$ has a left-inversion.
\end {Thm}

\begin {Def}
A {\em basis} $\{ e_{ \lambda} \setmid \lambda \in \Lambda \}$ of a $k$-module $M$ is a minimal system of generators (i.e. there is no $\lambda_0 \in \Lambda$ such that the elements $\{ e_{ \lambda} \setmid \lambda \in \Lambda \setminus \{ \lambda_0 \} \}$ generate $M$). A subset $S \subset M$ {\em generate} $M$ if any element of $M$ is written as a linear combination
\[
a_1 \dx v_1 \px \cdots \px a_r \dx v_r
\]
of elements of $S$ over $k$.
\end {Def}

\begin {Def}
An element $e \in M \setminus \{ \mx \}$ is {\em extremal} if for any $v_1, v_2 \in M$ such that $v_1 \px v_2= e$ it implies $v_1= e$ or $v_2= e$. $M$ is {\em extremally generated} if $M$ is generated by extremal elements. An {\em extremal ray} of $M$ is the submodule generated by an extremal element of $M$.
\end {Def}

\begin {Def}
The {\em dimension} of a straight reflexive $k$-module $M$ is the number of extremal rays. 
\end {Def}

The number of extremal rays of $M$ is equal to the number of elements of any basis of $M$. The proof is given in section \ref{Sec:FA6}.

\begin {Thm} \label {Prop:FE2}
Let $\alpha \colon M \rightarrow N$ be an injective homomorphism of finitely generated straight reflexive $k$-modules. Then
\begin {enumerate}
\renumi{\statementenumi}
\item $\xxdim( M) \leq \xxdim( N)$.
\item If $\xxdim( M)= \xxdim( N)$, then $\alpha$ is lightly surjective.
\end {enumerate}
\end {Thm}

\begin {Thm} \label {Prop:FE3}
Let $M$ be a finitely generated straight pre-reflexive $k$-module. Then $M$ is reflexive.
\end {Thm}

\begin {Thm} \label {Prop:FE4}
Let $\alpha \colon M \rightarrow N$ be an injective homomorphism of straight pre-reflexive $k$-modules. Suppose that $M$ has a basis, and that $N$ is finitely generated. Then $M$ is finitely generated.
\end {Thm}

Let $P$ be a polytope in $\rrt \pp^n$. $P$ is the tropically convex hull of finitely many points $p_1, \ldots, p_r$. Let
\[
\ph \colon \rrt^{ n+1} \setminus \{ \mx \} \rmap \rrt \pp^n
\]
be the canonical projection. Then the subset
\[
M= \inv{ \ph}( P) \cup \{ \mx \} \subset \rrt^{ n+1}
\]
is a submodule generated by elements $v_1, \ldots, v_r$ such that $\ph( v_i)= p_i$ ($1 \leq i \leq r$). Also we have an injection
\[
\iota \colon \rrt^n \rmap \rrt \pp^n
\]
defined by $(a_1, \ldots, a_n) \mapsto (0, a_1, \ldots, a_n)$. This map induces an embedding $\rr^n \subset \rrt^n \subset \rrt \pp^n$. A polytope $P \subset \rrt \pp^n$ is said to be a polytrope if it is a real convex subset of $\rr^n$.

\begin {Thm} \label {Prop:FE5}
Let $P$ be a polytope in $\rrt \pp^n$ with the corresponding submodule $M \subset \rrt^{ n+1}$.
\begin {enumerate}
\renumi{\statementenumi}
\item If $P$ is a polytrope, then $M$ is straight reflexive.
\item If $M$ is straight reflexive, then $P$ is the tropically convex hull of at most $n+1$ points.
\end {enumerate}
\end {Thm}

Let $C$ be a tropical curve. Let $D$ be a divisor on $C$. Let $H^0( C, \ooo_{ C}( D))$ be the set of the sections of $D$. (A section of $D$ is a rational function $f \colon C \rightarrow \rrt$ such that either $f= \mx$ or $( f)+ D \geq 0$.) For $r \in \zz_{ \geq 0}$, let
\[
U( D, r)= C^r \setminus S( D, r),
\]
\[
S( D, r)= \{( P_1, \ldots, P_r) \in C^r \setmid H^0( C, \ooo_{ C}( D- \sum_{ 1 \leq i \leq r} P_i)) \not= \mx \}.
\]
Let $U( D, r)= \emptyset$ if $r= -1$. The following theorem is known.

\begin {Thm}[Gathmann and Kerber \cite{Gath0612}] \label {Prop:FE6}
Let $C$ be a compact tropical curve with first Betti number $b_1( C)$. Let $D$ be a divisor on $C$. Let $K$ be the canonical divisor on $C$. Then
\[
r( D)- r( K- D)= 1- b_1( C)+ \deg(D),
\]
where
\[
r( D)= \max \{ r \in \zz_{ \geq -1} \setmid U( D, r) = \emptyset \}.
\]
\end {Thm}

The set $M= H^0( C, \ooo_{ C}( D))$ is a $\rrt$-module with addition
\[
(f \px g)( P)= f( P) \px g( P)
\]
and scalar multiplication
\[
(a \dx f)( P)= a \dx f( P).
\]
The dimension of $M$ is defined as follows.

\begin {Def}
The {\em dimension} of a $k$-module $M$ is the maximum dimension of the straight reflexive submodules of $M$.
\end {Def}

This definition is compatible with the previous one. If $M$ is straight reflexive, then the maximum dimension of the straight reflexive submodules of $M$ equals the dimension of $M$ by Theorem \ref{Prop:FE2}.

\begin {Thm} \label {Prop:FE7}
Let $C$ be a tropical curve. Let $D$ be a divisor on $C$. Then the inequality
\[
r( D) \leq \xxdim H^0(C, \ooo_{ C}( D))- 1
\]
is fulfilled.
\end {Thm}

\section {Tropical algebra} \label {Sec:FB2}

\subsection {Tropical semigroups, semirings, and semifields} \label {Sec:FA4}

\begin {Prop} \label {Prop:FA1}
Let $M$ be a tropical semigroup. Then there is a unique partial order relation `$\lqx$' such that for any $v, w \in M$ it implies
\[
\xxsup \{v, w \}= v \px w.
\]
\end {Prop}

\begin {proof}
We define a relation `$\lqx$' on $M$ as follows.
\[
v \lqx w \iff v \px w= w.
\]
This is a partial order relation, because $v \px v= v$. The element $v \px w$ is the minimum element of the upper bounds of $\{ v, w \}$.
\end {proof}

Let $A$ be a tropical semiring.

\begin {Ex}
The {\em semiring of polynomials} $B= A[x_1, \ldots, x_n]$ is the set of polynomials
\begin {align*}
f & = \bigpx_{ i} a_i \dx x^{ \dx i} \\
& = \bigpx_{i_1, \ldots, i_n \geq 0} a_{ i_1 \ldots i_n} \dx x_1^{ \dx i_1} \dx \cdots \dx x_n^{ \dx i_n}
\end {align*}
with coefficients $a_i \in A$, equipped with addition and multiplication of polynomials. $B$ is a tropical semiring. An element $f \in B$ is said to be a {\em tropical polynomial over} $A$. The induced map
\[
\begin {array}[c]{rcl}
\overline{ f} \colon A^n & \rmap & A \\
( a_1, \ldots, a_n) & \mapsto & f( a_1, \ldots, a_n)
\end {array}
\]
is said to be a tropical polynomial function.
\end {Ex}

\begin {Remark}
We use the notation $m a$ by the meaning of tropical $m$-th power $a^{ \dx m}$. For example, $2( a \px b)$ means the second power of $( a \px b)$, so we have
\begin {align*}
2( a \px b) & = 2 a \px a \dx b \px a \dx b \px 2 b \\
& = 2 a \px a \dx b \px 2 b.
\end {align*}
Also a tropical polynomial is written as
\[
f= \bigpx_{ i} a_i \dx i x.
\]
\end {Remark}

\begin {Prop} \label {Prop:FC1}
Let $A$ be a tropical semiring. Let $f \in A[x_1, \ldots, x_n]$. Then for any $v, w \in A^n$,
\[
f( v \px w) \gqx f( v) \px f( w).
\]
\end {Prop}

\begin {proof}
Assume that
\[
f = i_1 x_1 \dx \cdots \dx i_n x_n, \]\[
v =( a_1, \ldots, a_n), \]\[
w =( b_1, \ldots, b_n).
\]
Then
\begin {align*}
f( v \px w) & = i_1 ( a_1 \px b_1) \dx \cdots \dx i_n ( a_n \px b_n) \\
& \gqx ( i_1 a_1 \dx \cdots \dx i_n a_n) \px( i_1 b_1 \dx \cdots \dx i_n b_n) \\
& = f( v) \px f( w).
\end {align*}
\end {proof}

Let $k$ be a tropical semifield. Recall that $k$ is said to be rational if it satisfies the following conditions.
\begin {enumerate}
\renumi{\conditionenumi}
\item $a \in k$, $m \in \nn$ $\then$ $\setexists b \in k$, $a= b^{ \dx m}$.
\item $k$ has no maximum element.
\end {enumerate}

\begin {Prop} \label {Prop:FB1}
Let $k$ be a rational tropical semifield. Then for any $a \in k$ it implies
\[
\inf_{ k} \{ b \in k \setmid a \ltx b \}= a.
\]
\end {Prop}

\begin {proof}
The case of $a= \mx$. Suppose that there is an element $c \in k \setminus \{ \mx \}$ such that $\upperx{ k}( c)= k \setminus \{ \mx \}$. Then the element $0 \sx c$ is the maximum element of $k$, which is contradiction.
\par
The case of $a \not= \mx$. The condition $a \ltx b$ is fulfilled if and only if $0 \ltx b \sx a$. So we may assume $a= 0$. Suppose that there is an element $c \nlqx 0$ such that $c$ is a lower-bound of the set $\{ b \in k \setmid 0 \ltx b \}$. There is an element $c' \in k$ such that $c= (c')^{ \dx 2}= 2 c'$. Since $0 \ltx 0 \px c'$, we have $c \lqx 0 \px c'$. So we have
\begin {align*}
2 ( 0 \px c') & = 0 \px c' \px 2 c' \\
& = 0 \px c' \px c \\
& = 0 \px c', \\
0 \px c' & = 0.
\end {align*}
So we have $c \lqx 0$, which is contradiction.
\end {proof}

\subsection {Modules over a tropical semifield}

Let $k$ be a tropical semifield. Let $M$ be a $k$-module.

\begin {Def}
A {\em submodule} $N$ of $M$ is a subset with the following conditions.
\begin {enumerate}
\renumi{\conditionenumi}
\item $\mx \in N$.
\item If $v, w \in N$ then $v \px w \in N$.
\item If $v \in N$ and $a \in k$ then $a \dx v \in N$.
\end {enumerate}
\end {Def}

\begin {Ex} \label {Ex:FA3}
Suppose that $k$ is totally ordered. Let $q \in k[x_1, \ldots, x_n]$ be a homogeneous polynomial of degree $m$. Let $p \colon k^n \rightarrow k$ be a homomorphism of $k$-modules. Then the subset
\[
M= \{ v \in k^n \setmid m p( v) \lqx q( v) \}
\]
is a submodule of $k^n$. Indeed, for $v, w \in M$ and $a \in k$,
\begin {align*}
m p( a \dx v) & = m( a \dx p( v)) \\
& = m a \dx m p( v) \\
& \lqx m a \dx q( v) \\
& = q( a \dx v),
\end {align*}
\begin {align*}
m p( v \px w) & = m( p( v) \px p( w)) \\
& = \max \{ m p( v), m p( w) \} \\
& = m p( v) \px m p( w) \\
& \lqx q( v) \px q( w).
\end {align*}
By Proposition \ref{Prop:FC1}, we have $q( v) \px q( w) \lqx q( v \px w)$.
\end {Ex}

\begin {Ex}
A free module $M= k^n$ of finite rank is reflexive. Indeed there is a pairing map $\pairx{ \sethyphen}{ \sethyphen} \colon k^n \times k^n \rightarrow k$ defined by
\[
\pairx{ (a_1, \ldots, a_n)}{ (b_1, \ldots, b_n)}= a_1 \dx b_1 \px \cdots \px a_n \dx b_n.
\]
So we have $\che{( k^n)} \iso k^n$.
\end {Ex}

Recall that $M$ is said to be pre-reflexive if the homomorphism $\iota_{ M} \colon M \rightarrow \che{( \che{ M})}$ is injective.

\begin {Prop} \label {Prop:LK8}
$M$ is pre-reflexive if and only if there is an injection $M \rightarrow F$ for some direct product $F= \bigdirdx_{ \lambda \in \Lambda} k$.
\end {Prop}

\begin {proof}
There is an injection $\che{( \che{ M})} \rightarrow \bigdirdx_{ \lambda \in \Lambda} k$, where $\Lambda$ is the set $\che{ M}$. Conversely, if there is an injection $M \rightarrow F$ for some direct product $F$, then $M$ is pre-reflexive, because $F$ is pre-reflexive.
\end {proof}

\begin {Lem} \label {Prop:FB3}
Suppose that $k$ is rational. Let $M$ be a pre-reflexive $k$-module. Then for any $v \in M$ and any $a \in k$ it implies
\[
\inf_{ M} \{ b \dx v \setmid b \in k, a \ltx b \}= a \dx v.
\]
\end {Lem}

\begin {proof}
Let $w \in M$ be a lower-bound of the subset $\{ b \dx v \setmid b \in k, a \ltx b \}$. For $\xi \in \che{ M}$ and $b \in k$ such that $a \ltx b$, we have
\[
\xi( w) \lqx b \dx \xi( v).
\]
By Proposition \ref{Prop:FB1}, we have
\[
\xi( w) \lqx a \dx \xi( v).
\]
Since $M$ is pre-reflexive, we have $w \lqx a \dx v$.
\end {proof}

\begin {Lem} \label {Prop:LD11}
Suppose that $k$ is totally ordered. Let $M$ be a pre-reflexive $k$-module. Then for any $v, w \in M$ and any $a \in k$,
\[
v \nlqx w, a \ltx 0 \then v \nlqx w \px a \dx v.
\]
\end {Lem}

\begin {proof}
Since $M$ is pre-reflexive, there is an element $\xi \in \che{ M}$ such that $\xi( v) \nlqx \xi( w)$. Since $k$ is totally ordered, we have $\xi( w) \ltx \xi( v)$. So
\[
\max \{ \xi( w), a \dx \xi( v) \} \ltx \xi( v).
\]
So we have the conclusion.
\end {proof}

\begin {Ex} \label {Ex:LA3}
Let $G$ be a tropical semigroup with at least two elements. Let $M=( G \times \rr) \cup \{ \mx \}$ be the $\rrt$-module with addition
\[
( v, a) \px ( w, b)=( v \px w, a \px b)
\]
and scalar multiplication
\[
c \dx( v, a)= \begin {cases}
( v, c \dx a) & \text{if $c \in \rr$} \\
\mx & \text{if $c= \mx$}.
\end {cases}
\]
$M$ is a $\rrt$-module generated by the subset $G \times \{ 0 \}$. $M$ is not pre-reflexive, because it does not satisfy Lemma \ref{Prop:LD11}. Let $v, w \in G$ be elements such that $v \nlqx w$. Then
\begin {align*}
( v, 0) & \nlqx ( w, 0), \\
( v, 0) & \lqx ( w, 0) \px( -1) \dx( v, 0).
\end {align*}
\end {Ex}

\subsection {Basis and extremal rays} \label {Sec:FA6}

Let $k$ be a totally ordered tropical semifield. Let $M$ be a $k$-module. Recall that an element $e \in M \setminus \{ \mx \}$ is said to be extremal if for any $v_1, v_2 \in M$ such that $v_1 \px v_2= e$ it implies $v_1= e$ or $v_2= e$.

\begin {Prop} \label {Prop:LK6}
Let $M$ be a pre-reflexive $k$-module. Then the following are equivalent.
\begin {enumerate}
\renumi{\conditionenumi}
\item There is a basis of $M$.
\item $M$ is extremally generated.
\end {enumerate}
More precisely, a system of generators $E= \{ e_{ \lambda} \setmid \lambda \in \Lambda \}$ is a basis if and only if each $e_{ \lambda}$ is extremal and it satisfies $k \dx e_{ \lambda} \not= k \dx e_{ \mu}$ ($\lambda \not= \mu$).
\end {Prop}

\begin {proof}
Suppose that there is a basis $E$ of $M$. Let $e_1$ be an element of the basis $E$. Let $v_1, v_2 \in M$ be elements such that $v_1 \px v_2= e_1$. There are elements $e_2, e_3, \ldots, e_r$ of the basis $E$ and elements $a_i, b_i \in k$ such that
\begin {align*}
v_1 & = a_1 \dx e_1 \px a_2 \dx e_2 \px \cdots \px a_r \dx e_r, \\
v_2 & = b_1 \dx e_1 \px b_2 \dx e_2 \px \cdots \px b_r \dx e_r.
\end {align*}
Since $k$ is totally ordered, we may assume $a_1 \lqx b_1$.
Then
\[
e_1= b_1 \dx e_1 \px w,
\]
where
\[
w=( a_2 \px b_2) \dx e_2 \px \cdots  \px( a_r \px b_r) \dx e_r.
\]
Since $E$ is a basis, we have $w \not= e_1$. By Lemma \ref{Prop:LD11}, we have $b_1= 0$. It means $v_2 \gqx e_1$. So we have $v_2= e_1$. Thus $e_1$ is extremal.
\par
Conversely, let $E$ be a system of generators that consists of extremal elements with different extremal rays. Suppose that $E$ is not a basis. There are elements $e_1, e_2, \ldots, e_r$ of $E$ and elements $a_i \in k$ such that
\[
e_1= a_2 \dx e_2 \px \cdots \px a_r \dx e_r.
\]
Since $e_1$ is extremal, there is a number $i$ such that $e_1= a_i \dx e_i$, which is contradiction.
\end {proof}

\begin {Prop} \label {Prop:FA4}
Let $\alpha \colon M \rightarrow N$ be a homomorphism of $k$-modules. Let $w \in N$ be an extremal element. Then any minimal element of the subset $\inv{ \alpha}( w)$ is extremal.
\end {Prop}

\begin {proof}
Let $e \in M$ be a minimal element of $\inv{ \alpha}( w)$. Let $v_1, v_2 \in M$ be elements such that $v_1 \px v_2= e$. Then $\alpha( v_1) \px \alpha( v_2)= w$. Since $w$ is extremal, we may assume $\alpha( v_1)= w$. Then $v_1$ is a lower-bound of $e$ in $\inv{ \alpha}( w)$. Since $e$ is minimal, we have $v_1= e$.
\end {proof}

\subsection {Locators}

Let $k$ be a totally ordered tropical semifield. Let $M$ be a $k$-module. For a subset $S \subset M$, the lower-saturation $\lowerx{ M}( S)$ is defined by
\[
\lowerx{ M}( S)= \bigcup_{ w \in S} \{ v \in M \setmid v \lqx w \}.
\]
The set of the lower-bounds $\lowerbx{ M}( S)$ is defined by
\[
\lowerbx{ M}( S)= \bigcap_{ w \in S} \{ v \in M \setmid v \lqx w \}.
\]
A subset $S \subset M$ is said to be lower-saturated if $\lowerx{ M}( S) = S$.

\begin {Def}
A {\em locator} $S$ of $M$ is a lower-saturated subsemigroup of the semigroup $( M, \px)$ that generates the $k$-module $M$.
\end {Def}

Let $\cheloc{ M}$ denote the set of the locators of a $k$-module $M$, equipped with addition
\[
S \chepx T= S \cap T
\]
and scalar multiplication
\[
a \chedx S= \begin {cases}
(\sx a) \dx S & \text{if $a \in k \setminus \{ \mx \}$} \\
M & \text{if $a= \mx$}.
\end {cases}
\]

\begin {Prop} \label {Prop:LC2}
$( \cheloc{ M}, \chepx, \chedx)$ is a $k$-module with zero element $M$. There is a homomorphism
\[
i \colon \che{ M} \rmap \cheloc{ M}
\]
defined by
\[
i( \xi)= \{ v \in M \setmid \pairx{ v}{ \xi} \lqx 0 \}.
\]
\end {Prop}

\begin {proof}
$\cheloc{ M}$ is a tropical semigroup. Indeed,
\[
S \chepx S= S \cap S= S.
\]
$\cheloc{ M}$ is a $k$-module. Indeed, for $a, b \in k$ such that $a \lqx b$, since $S$ is lower-saturated, we have
\[
\sx b \dx S \subset \sx a \dx S.
\]
So we have
\begin {align*}
(a \px b) \chedx S & = \sx b \dx S \\
& =( \sx a \dx S) \cap( \sx b \dx S) \\
& = a \chedx S \px b \chedx S.
\end {align*}
$i$ is a homomorphism. Indeed, for $v \in M$,
\begin {align*}
v \in i( \xi_1 \px \xi_2) & \iff \pairx{ v}{ \xi_1 \px \xi_2} \lqx 0 \\
& \iff \pairx{ v}{ \xi_1} \px \pairx{ v}{ \xi_2} \lqx 0 \\
& \iff v \in i( \xi_1) \cap i( \xi_2) \\
& \iff v \in i( \xi_1) \chepx i( \xi_2).
\end {align*}
So $i( \xi_1 \px \xi_2)= i( \xi_1) \chepx i( \xi_2)$.
\begin {align*}
v \in i( a \dx \xi) & \iff \pairx{ v}{ a \dx \xi} \lqx 0 \\
& \iff \pairx{ a \dx v}{ \xi} \lqx 0 \\
& \iff a \dx v \in i( \xi) \\
& \iff v \in a \chedx i( \xi).
\end {align*}
So $i( a \dx \xi)= a \chedx i( \xi)$.
\end {proof}

\begin {Lem} \label {Prop:FB2}
Suppose that $k$ is quasi-complete and rational.
\begin {enumerate}
\renumi{\statementenumi}
\item For any locator $S \in \cheloc{ M}$ there is a unique element $\xi \in \che{ M}$ that satisfies the following conditions.
\begin {align*}
\pairx{ v}{ \xi} & \lqx 0 \quad (v \in S), \\
\pairx{ v}{ \xi} & \gqx 0 \quad (v \in M \setminus S).
\end {align*}
\item The mapping $S \mapsto \xi$ induces a homomorphism
\[
p \colon \cheloc{ M} \rmap \che{ M}
\]
which satisfies $p \circ i= \xxid_{ \che{ M}}$.
\end {enumerate}
\end {Lem}

\begin {proof}
\statement{1} Let $\xi \colon M \rightarrow k$ be the map defined as follows.
\[
\xi( v)= \inf_{ k} \{ a \in k \setmid v \in a \dx S \}.
\]
The set in right side is non-empty. (Since $S$ generates the $k$-module $M$, there are $s_i \in S$ and $a_i \in k$ such that
\[
v= a_1 \dx s_1 \px \cdots \px a_r \dx s_r.
\]
Let $a$ be the maximum element of $a_1, \ldots, a_r$. Since $S$ is lower-saturated, there are $s'_i \in S$ such that
\[
v= a \dx( s'_1 \px \cdots \px s'_r).
\]
Since $S$ is a subsemigroup, we have $v \in a \dx S$.) For any $v \in M \setminus S$ we have $\xi( v) \gqx 0$, because $S$ is lower-saturated. For any $v \in S$, we have $\xi( v) \lqx 0$.
\par
We show that $\xi$ is a homomorphism. Since $S$ is lower-saturated, we have
\[
\xi( v) \px \xi( w) \lqx \xi( v \px w).
\]
Suppose that $\xi( v) \px \xi( w) \ltx \xi( v \px w)$. There are $a, b \in k$ such that $a \px b \ltx \xi( v \px w)$ and $v \in a \dx S$ and $w \in b \dx S$. Then $v \px w \in( a \px b) \dx S$. So we have $\xi( v \px w) \lqx a \px b$, which is contradiction.
\par
We prove uniqueness. Let $\xi \in \che{ M}$ be an element that satisfies the following conditions.
\begin {align*}
\pairx{ v}{ \xi} & \lqx 0 \quad (v \in S), \\
\pairx{ v}{ \xi} & \gqx 0 \quad (v \in M \setminus S).
\end {align*}
Then
\begin {align*}
\pairx{ v}{ \xi} & \lqx \inf_{ k} \{ a \in k \setmid v \in a \dx S \} \\
& \lqx \inf_{k} \{ a \in k \setmid \pairx{ v}{ \xi} \ltx a \}.
\end {align*}
By Proposition \ref{Prop:FB1},
\[
\inf_{k} \{ a \in k \setmid \pairx{ v}{ \xi} \ltx a \}= \pairx{ v}{ \xi}.
\]
So we have
\[
\pairx{ v}{ \xi}= \inf_{ k} \{ a \in k \setmid v \in a \dx S \}.
\]
\par
\statement{2} We have
\begin {align*}
\pairx{ v}{ p( S) \px p( T)} & \lqx 0 \quad (v \in S \cap T), \\
\pairx{ v}{ p( S) \px p( T)} & \gqx 0 \quad (v \in M \setminus( S \cap T)).
\end {align*}
It means $p( S) \px p( T)= p( S \chepx T)$. So $p$ is a homomorphism. For $\xi \in \che{ M}$, let
\[
S= \{ v \in M \setmid \pairx{ v}{ \xi} \lqx 0 \}.
\]
Then
\begin {align*}
\pairx{ v}{ \xi} & \lqx 0 \quad (v \in S), \\
\pairx{ v}{ \xi} & \gqx 0 \quad (v \in M \setminus S).
\end {align*}
It means $\xi= p( S)$.
\end {proof}

\subsection {Straight modules}

Let $k$ be a totally ordered tropical semifield. Recall that a $k$-module $M$ is said to be straight if it satisfies the following conditions.
\begin {enumerate}
\renumi{\conditionenumi}
\item Any two elements $v, w \in M$ admit the infimum $\infx{ M}{ v}{ w}$.
\item $v_1, v_2, w \in M \then \infx{ M}{ v_1 \px v_2}{ w}= \infx{ M}{ v_1}{ w} \px \infx{ M}{ v_2}{ w}$.
\item $v_1, v_2, w \in M \then \infx{ M}{ v_1}{ v_2} \px w= \infx{ M}{ v_1 \px w}{ v_2 \px w}$.
\end {enumerate}

\begin {Prop} \label {Prop:LI3}
The above conditions \condition{ii}, \condition{iii} are equivalent.
\end {Prop}

\begin {proof}
\condition{ii} $\then$ \condition{iii}.
\begin {align*}
\infx{ M}{ v_1 \px w}{ v_2 \px w} & = \infx{ M}{ v_1}{ v_2} \px \infx{ M}{ v_1}{ w} \px \infx{ M}{ w}{ v_2} \px w \\
& = \infx{ M}{ v_1}{ v_2} \px w.
\end {align*}
\condition{iii} $\then$ \condition{ii} is similar.
\end {proof}

\begin {Def}
A homomorphism $\alpha \colon M \rightarrow N$ of $k$-modules is {\em lattice-preserving} if for any $v, w \in M$ and any lower-bound $x \in \lowerbx{ N}( \alpha( v), \alpha( w))$ there is a lower-bound $y \in \lowerbx{ M}( v, w)$ such that $x \lqx \alpha( y)$.
\end {Def}

If $M$, $N$ are ordered lattices, $\alpha$ is lattice-preserving if and only if it preserves the infimum of any two elements.

\begin {Prop} \label {Prop:FA5}
Let $\alpha \colon M \rightarrow N$ be a lattice-preserving injective homomorphism of $k$-modules such that $N$ is straight. Then $M$ is straight.
\end {Prop}

\begin {proof}
For $v, w \in M$, let $x= \infx{ N}{ \alpha( v)}{ \alpha( w)}$. There is a lower-bound $y$ of $\{ v, w \}$ such that $x \lqx \alpha( y)$. Then $y= \infx{ M}{ v}{ w}$. (Let $y' \in M$ be a lower-bound of $\{ v, w \}$. Then $\alpha( y') \lqx x \lqx \alpha( y)$. Since $\alpha$ is injective, we have $y' \lqx y$.) $\alpha( y)$ is a lower-bound of $\{ \alpha( v), \alpha( w) \}$. So we have $x= \alpha( y)$. $M$ is finitely distributive, because $\alpha$ preserves the infimum of any two elements.
\end {proof}

\begin {Prop} \label {Prop:FA3}
Suppose that $k$ is quasi-complete and rational. Let $M$ be a straight $k$-module. Then $\che{ M}$ and $\cheloc{ M}$ are straight.
\end {Prop}

\begin {proof}
We show that $\cheloc{ M}$ is straight. For $S, T \in \cheloc{ M}$, let
\[
U= S \px T= \{ s \px t \setmid s \in S, t \in T \}.
\]
$U$ is lower-saturated. (Let $v \in M$ and $s \in S$ and $t \in T$ be elements such that $v \lqx s \px t$. Then
\begin {align*}
v & = \infx{ M}{ v}{ s \px t} \\
& = \infx{ M}{ v}{ s} \px \infx{ M}{ v}{ t}.
\end {align*}
So we have $v \in U$.) $U$ is a locator of $M$, and we have
\[
U= \infx{ \cheloc{ M}}{ S}{ T}.
\]
$\cheloc{ M}$ is finitely distributive. Indeed,
\[
( S_1 \cap S_2) \px T=( S_1 \px T) \cap( S_2 \px T).
\]
(Let $v$ be an element of right side. There are $s_1 \in S_1$ and $s_2 \in S_2$ and $t_1, t_2 \in T$ such that
\[
v= s_1 \px t_1= s_2 \px t_2.
\]
Then 
\begin {align*}
v & = \infx{ M}{ s_1 \px t_1}{ s_2 \px t_2} \\
& = \infx{ M}{ s_1}{ s_2} \px \infx{ M}{ s_1}{ t_2} \px \infx{ M}{ t_1}{ s_2} \px \infx{ M}{ t_1}{ t_2}.
\end {align*}
So we have $v \in ( S_1 \cap S_2) \px T$.)
\par
We show that $\che{ M}$ is straight. For $\xi_1, \xi_2 \in \che{ M}$, let $S_1, S_2 \in \cheloc{ M}$ be the induced element. There is a unique element $\eta \in \che{ M}$ that satisfies the following conditions (Lemma \ref{Prop:FB2}).
\begin {align*}
\pairx{ v}{ \eta} & \lqx 0 \quad (v \in S_1 \px S_2), \\
\pairx{ v}{ \eta} & \gqx 0 \quad (v \in M \setminus( S_1 \px S_2)).
\end {align*}
We have $\eta= \infx{ \che{ M}}{ \xi_1}{ \xi_2}$. So the canonical injection $i \colon \che{ M} \rightarrow \cheloc{ M}$ is lattice-preserving. Since $\cheloc{ M}$ is straight, $\che{ M}$ is straight (Proposition \ref{Prop:FA5}).
\end {proof}

\begin {Prop} \label {Prop:FA6}
Let $M$ be a $k$-module. Let $\eta \colon M \rightarrow k$ be a lattice-preserving homomorphism. Then $\eta$ is an extremal element of $\che{ M}$.
\end {Prop}

\begin {proof}
Suppose that $\eta$ is not extremal. There are elements $\xi_1, \xi_2 \in \che{ M}$ and elements $v_1, v_2 \in M$ such that $\xi_1 \px \xi_2= \eta$ and $\pairx{ v_1}{ \xi_1} \ltx \pairx{ v_1}{ \eta}$ and $\pairx{ v_2}{ \xi_2} \ltx \pairx{ v_2}{ \eta}$. We may assume $\pairx{ v_1}{ \eta}= \pairx{ v_2}{ \eta}= 0$. Since $\eta$ is lattice-preserving, there is a lower-bound $w$ of $\{ v_1, v_2 \}$ such that $\pairx{ w}{ \eta}= 0$. Then
\begin {align*}
0 & = \pairx{ w}{ \eta} \\
& = \pairx{ w}{ \xi_1 \px \xi_2} \\
& \lqx \pairx{ v_1}{ \xi_1} \px \pairx{ v_2}{ \xi_2} \\
& \ltx \pairx{ v_1}{ \eta} \px \pairx{ v_2}{ \eta} \\
& = 0,
\end {align*}
which is contradiction.
\end {proof}

\begin {Def}
A {\em dual element} $\eta \in \che{ M}$ of an element $e \in M$ is an element with the following conditions.
\begin {enumerate}
\renumi{\conditionenumi}
\item $\pairx{ e}{ \eta} = 0$.
\item $v \in M, \xi \in \che{ M} \then \pairx{ v}{ \eta} \dx \pairx{ e}{ \xi} \lqx \pairx{ v}{ \xi}$.
\end {enumerate}
\end {Def}

\begin {Prop} \label {Prop:FA7}
The dual element of an element $e \in M$ is unique.
\end {Prop}

\begin {proof}
Let $\eta$ be a dual element of $e$. Then
\[
\eta= \min \{ \xi \in \che{ M} \setmid \pairx{ e}{ \xi}= 0 \},
\]
because
\[
\eta \dx \pairx{ e}{ \xi} \lqx \xi
\]
for any $\xi \in \che{ M}$.
\end {proof}

\begin {Prop} \label {Prop:FA8}
Let $e \in M$ be an element of a pre-reflexive $k$-module $M$. Suppose that $e$ has the dual element $\eta \in \che{ M}$. Then
\begin {enumerate}
\renumi{\statementenumi}
\item $e$ is an extremal element.
\item $\eta \colon M \rightarrow k$ is a lattice-preserving homomorphism (therefore is an extremal element of $\che{ M}$).
\end {enumerate}
\end {Prop}

\begin {proof}
\statement{1} Let $v_1, v_2 \in M$ be elements such that $v_1 \px v_2= e$. Then
\[
\pairx{ v_1}{ \eta} \px \pairx{ v_2}{ \eta}= \pairx{ e}{ \eta}= 0.
\]
We may assume $\pairx{ v_1}{ \eta}= 0$. For $\xi \in \che{ M}$,
\begin {align*}
\pairx{ e}{ \xi} & = \pairx{ v_1}{ \eta} \dx \pairx{ e}{ \xi} \\
& \lqx \pairx{ v_1}{ \xi}.
\end {align*}
Since $M$ is pre-reflexive, we have $e \lqx v_1$. So we have $e= v_1$.
\par
\statement{2} Let $v_1, v_2 \in M$ be elements such that $\pairx{ v_1}{ \eta} \lqx \pairx{ v_2}{ \eta}$. Let $w= \pairx{ v_1}{ \eta} \dx e$. For $\xi \in \che{ M}$, we have $\pairx{ w}{ \xi} \lqx \pairx{ v_i}{ \xi}$. Since $M$ is pre-reflexive, we have $w \lqx v_i$. So $w$ is a lower-bound of $\{ v_1, v_2 \}$ such that $\pairx{ w}{ \eta}= \pairx{ v_1}{ \eta}$. Thus $\eta$ is lattice-preserving. By Proposition \ref{Prop:FA6}, $\eta$ is an extremal element of $\che{ M}$.
\end {proof}

\begin {Lem} \label {Prop:FA9}
Suppose that $k$ is quasi-complete and rational. Let $M$ be a straight pre-reflexive $k$-module. Then any extremal element of $M$ has the dual element.
\end {Lem}

\begin {proof}
Let $e \in M$ be an extremal element. The subset
\[
S= \{ v \in M \setmid e \nlqx v \}
\]
is a subsemigroup. (Let $v_1, v_2 \in M$ be elements such that $e \lqx v_1 \px v_2$. Then
\begin {align*}
e & = \infx{ M}{ e}{ v_1 \px v_2} \\
& = \infx{ M}{ e}{ v_1} \px \infx{ M}{ e}{ v_2}.
\end {align*}
We may assume $e= \infx{ M}{ e}{ v_1}$. Then $e \lqx v_1$.) Also $S$ generates the $k$ module $M$. (Let $v \in M$ be any element. By Lemma \ref{Prop:FB3}, we have
\[
\inf_{ M} \{ b \dx v \setmid b \in k \setminus \{ \mx \} \}= \mx.
\]
So there is $b \in k \setminus \{ \mx \}$ such that $e \nlqx b \dx v$.) Thus $S$ is a locator of $M$. By Lemma \ref{Prop:FB2}, there is a unique element $\eta \in \che{ M}$ that satisfies the following conditions.
\begin {align*}
\pairx{ v}{ \eta} & \lqx 0 \quad (v \in S), \\
\pairx{ v}{ \eta} & \gqx 0 \quad (v \in M \setminus S).
\end {align*}
Then
\begin {align*}
\pairx{ v}{ \eta} & \lqx \inf_{ k} \{ a \in k \setmid a \dx e \nlqx v \} \\
& \lqx \inf_{ k} \{ a \in k \setmid \pairx{ v}{ \eta} \ltx a \} \\
& = \pairx{ v}{ \eta}.
\end {align*}
So we have
\[
\pairx{ v}{ \eta}= \inf_{ k} \{ a \in k \setmid a \dx e \nlqx v \}.
\]
So
\begin {align*}
\pairx{ e}{ \eta} & = \inf_{ k} \{ a \in k \setmid 0 \ltx a \} \\
& = 0.
\end {align*}
Also, for any $a \in k$ such that $0 \ltx a$, we have
\[
( \pairx{ v}{ \eta} \sx a) \dx e \lqx v.
\]
By Lemma \ref{Prop:FB3}, we have
\[
\pairx{ v}{ \eta} \dx e \lqx v.
\]
Thus $\eta$ is the dual element of $e$.
\end {proof}

\begin {Lem} \label {Prop:FA10}
Let $M$ be a finitely generated pre-reflexive $k$-module. Let $\beta \colon k^n \rightarrow M$ be the surjection defined by a basis $\{ e_1, \ldots, e_n \}$ of $M$. Suppose that $e_i$ has the dual element $\eta_i$ ($1 \leq i \leq n$). Then
\begin {enumerate}
\renumi{\statementenumi}
\item $M$ is straight.
\item The homomorphism $\alpha \colon M \rightarrow k^n$ defined by the elements $\eta_1, \ldots, \eta_n$ is a right-inversion of $\beta$, i.e. $\beta \circ \alpha= \xxid_{ M}$.
\item $\alpha$ is the unique right-inversion of $\beta$.
\end {enumerate}
\end {Lem}

\begin {proof}
For $v \in M$, we have
\[
v \gqx \pairx{ v}{ \eta_1} \dx e_1 \px \cdots \px \pairx{ v}{ \eta_n} \dx e_n.
\]
It means $\beta \circ \alpha \lqx \xxid_{ M}$. Also, for $1 \leq i \leq n$, we have
\[
\beta \circ \alpha( e_i) \gqx \pairx{ e_i}{ \eta_i} \dx e_i= e_i.
\]
Since $M$ is generated by $\{ e_1, \ldots, e_n \}$, we have $\beta \circ \alpha= \xxid_{ M}$. So $\alpha$ is injective. Also $\alpha$ is lattice-preserving (Proposition \ref{Prop:FA8}). Since $k^n$ is straight, $M$ is straight (Proposition \ref{Prop:FA5}).
\par
We prove uniqueness. Let $\eta'_1, \ldots, \eta'_n \in \che{ M}$ be elements such that the induced homomorphism $M \rightarrow k^n$ is a right-inversion of $\beta$. Then we have
\[
v= \pairx{ v}{ \eta'_1} \dx e_1 \px \cdots \px \pairx{ v}{ \eta'_n} \dx e_n.
\]
So
\[
e_i = \pairx{ e_i}{ \eta'_i} \dx e_i \px w_i,
\]
where
\[
w_i= \bigpx_{ j \not= i} \pairx{ e_i}{ \eta'_j} \dx e_j.
\]
Since $\{ e_1, \ldots, e_n \}$ is a basis, we have $w_i \not= e_i$ and
\[
e_i = \pairx{ e_i}{ \eta'_i} \dx e_i
\]
(Proposition \ref{Prop:LK6}). So we have $\pairx{ e_i}{ \eta'_i}= 0$. Thus $\eta'_i$ is the dual element of $e_i$.
\end {proof}

\subsection {Existence of inversions}

Let $k$ be a totally ordered tropical semifield. Let $\alpha \colon M \rightarrow N$ be a homomorphism of $k$-modules.

\begin {Def}
An element $\xi \in \che{ M}$ {\em dominates} an element $w \in N$ if there is an element $v \in M$ such that $\pairx{ v}{ \xi} \lqx 0$ and $w \lqx \alpha( v)$.
\end {Def}

\begin {Prop} \label {Prop:FB4}
Let $\xi_i \in \che{ M}$ be an element that dominates $w_i \in N$ ($i= 1, 2$). Then any lower-bound $\xi \in \lowerbx{ \che{ M}}( \xi_1, \xi_2)$ dominates $w_1 \px w_2$.
\end {Prop}

\begin {proof}
There are elements $v_1, v_2 \in M$ such that $\pairx{ v_i}{ \xi_i} \lqx 0$ and $w_i \lqx \alpha( v_i)$. Then
\begin {align*}
\pairx{ v_1 \px v_2}{ \xi} & \lqx \pairx{ v_1}{ \xi_1} \px \pairx{ v_2}{ \xi_2} \\
& \lqx 0.
\end {align*}
Also we have $w_1 \px w_2 \lqx \alpha( v_1 \px v_2)$.
\end {proof}

Recall that a homomorphism $\alpha \colon M \rightarrow N$ is said to be lightly surjective if for any $w \in N$ there is $v \in M$ such that $w \lqx \alpha( v)$.

\begin {Lem} \label {Prop:LC14}
Let $\alpha \colon M \rightarrow N$ be an injective lightly surjective homomorphism of $k$-modules. Suppose that $\che{ M}$ is straight.
\begin {enumerate}
\renumi{\statementenumi}
\item There is a homomorphism $\gamma \colon N \rightarrow \cheloc{ \che{ M}}$ that satisfies the following condition. For any $w \in N$ the locator $\gamma( w)$ is the subsemigroup of $\che{ M}$ generated by the elements that dominates the element $w$.
\item The diagram
\[
\squarerd{ M}{ \alpha}{}{ N}{ \gamma}{ \che{( \che{ M})}}{ i}{ \cheloc{ \che{ M}}}
\]
commutes, i.e. for any $v \in M$ and any $\xi \in \che{ M}$ the condition $\pairx{ v}{ \xi} \lqx 0$ is fulfilled if and only if $\xi \in \gamma( \alpha( v))$.
\end {enumerate}
\end {Lem}

\begin {proof}
\statement{1}
For $w \in N$, let $\gamma( w) \subset \che{ M}$ be the subsemigroup of $\che{ M}$ generated by the elements that dominates the element $w$. $\gamma( w)$ is lower-saturated. (Let $\xi \in M$ and $\xi' \in \gamma( w)$ be elements such that $\xi \lqx \xi'$. There are elements $\xi_1, \ldots, \xi_r \in \che{ M}$ such that $\xi_i$ dominates $w$ and
\[
\xi'= \xi_1 \px \cdots \px \xi_r.
\]
Then
\begin {align*}
\xi & = \infx{ \che{ M}}{ \xi}{ \xi_1 \px \cdots \px \xi_r} \\
& = \infx{ \che{ M}}{ \xi}{ \xi_1} \px \cdots \px \infx{ \che{ M}}{ \xi}{ \xi_r}.
\end {align*}
So $\xi \in \gamma( w)$.) Also $\gamma( w)$ generates the $k$-module $\che{ M}$. (Let $\xi \in \che{ M}$ be any element. Since $\alpha$ is lightly surjective, there is $v \in M$ such that $w \lqx \alpha( v)$. Let $a \in k \setminus \{ \mx \}$ be an element such that $\pairx{ v}{ \xi} \lqx a$. Then $\sx a \dx \xi$ dominates $w$.) So $\gamma( w)$ is a locator of $\che{ M}$.
\par
We show that $\gamma$ is a homomorphism. For $w_1, w_2 \in N$, we have
\[
\gamma( w_1 \px w_2) \subset \gamma( w_1) \cap \gamma( w_2).
\]
Let $\xi$ be an element of right side. There are elements $\xi_{ i, j} \in \che{ M}$ ($1 \leq i \leq 2$, $1 \leq j \leq r$) such that $\xi_{i, j}$ dominates $w_i$ and
\[
\xi= \xi_{ 1, 1} \px \cdots \px \xi_{ 1, r}= \xi_{ 2, 1} \px \cdots \px \xi_{ 2, r}.
\]
Then
\begin {align*}
\xi & = \infx{ \che{ M}}{ \xi_{ 1, 1} \px \cdots \px \xi_{ 1, r}}{ \xi_{ 2, 1} \px \cdots \px \xi_{ 2, r}} \\
& = \bigpx_{ i, j} \eta_{ i, j},
\end {align*}
where
\[
\eta_{ i, j}= \infx{ \che{ M}}{ \xi_{ 1, i}}{ \xi_{ 2, j}}.
\]
$\eta_{ i, j}$ dominates $w_1 \px w_2$ (Proposition \ref{Prop:FB4}). So we have $\xi \in \gamma( w_1 \px w_2)$.
\par
\statement{2} Let $\xi \in \che{ M}$ be an element that dominates $\alpha( v)$. There is an element $v' \in M$ such that $\pairx{ v'}{ \xi} \lqx 0$ and $\alpha( v) \lqx \alpha( v')$. Since $\alpha$ is injective, we have $v \lqx v'$. So we have
\[
\pairx{ v}{ \xi} \lqx \pairx{ v'}{ \xi} \lqx 0.
\]
Let
\[
T= \{ \xi \in \che{ M} \setmid \pairx{ v}{ \xi} \lqx 0 \}.
\]
Now we have $\xi \in T$. Since $T$ is a subsemigroup, we have $\gamma( \alpha( v))= T$.
\end {proof}

\subsection {Straight reflexive modules} \label {Sec:FB1}

Let $k$ be a quasi-complete totally ordered rational tropical semifield. Recall that the dimension of a straight reflexive $k$-module $M$ is the number of extremal rays. By Proposition \ref{Prop:LK6}, the number of elements of any basis of $M$ is $\xxdim( M)$.

\begin {proof}[Proof of Theorem \ref{Prop:FE1}]
We have an isomorphism $\iota_{ M} \colon M \rightarrow \che{( \che{ M})}$ and a homomorphism $\gamma \colon N \rightarrow \cheloc{ \che{ M}}$ defined in Lemma \ref{Prop:LC14}. There is a left-inversion $p$ of the homomorphism $i \colon \che{( \che{ M})} \rightarrow \cheloc{ \che{ M}}$ (Lemma \ref{Prop:FB2}). By the commutative diagram
\[
\squarerd{ M}{ \alpha}{ \iota_{ M}}{ N}{ \gamma}{ \che{( \che{ M})}}{ i}{ \cheloc{ \che{ M}}}
\]
we have $\inv{ \iota_{ M}} \circ p \circ \gamma \circ \alpha= \xxid_{ M}$.
\end {proof}

\begin {proof}[Proof of Theorem \ref{Prop:FE2}]
By Lemma \ref{Prop:FA9} and Lemma \ref{Prop:FA10}, there is an injection $N \rightarrow k^n$, where $n= \xxdim( N)$. Let $N'$ be the lower-saturation of the image of $M \rightarrow k^n$. $N'$ is a free module of rank $n' \leq n$. If $n'= n$, then $\alpha$ is lightly surjective. Now we may assume that $N= k^n$ and that $\alpha$ is lightly surjective. By Theorem \ref{Prop:FE1}, $\alpha$ has a left-inversion $\beta \colon N \rightarrow M$. Since $\beta$ is surjective, we have $\xxdim( M) \leq \xxdim( N)$.
\end {proof}

\begin {proof}[Proof of Theorem \ref{Prop:FE3}]
By Lemma \ref{Prop:FA9} and Lemma \ref{Prop:FA10}, there is a right-inversion $\alpha \colon M \rightarrow k^n$ of the surjection $\beta \colon k^n \rightarrow M$. By the commutative diagram
\[
\xymatrix{
k^n \aru{r}{ \beta} \arisod & M \aru{d}{ \iota_{ M}} \\
k^n \ard{r}{ \che{( \che{ \beta})}} & \che{( \che{ M})}
}
\]
we have $\inv{ \iota_{ M}}= \beta \circ \che{( \che{ \alpha})}$.
\end {proof}

\begin {proof}[Proof of Theorem \ref{Prop:FE4}]
By Theorem \ref{Prop:FE3}, $N$ is reflexive. Similarly to the proof of Theorem \ref{Prop:FE2}, we may assume that $N= k^n$ and that $\alpha$ is lightly surjective. We have a homomorphism $\gamma \colon N \rightarrow \cheloc{ \che{ M}}$ defined in Lemma \ref{Prop:LC14}. There is a left-inversion $p$ of the homomorphism $i \colon \che{( \che{ M})} \rightarrow \cheloc{ \che{ M}}$ (Lemma \ref{Prop:FB2}). There is a homomorphism $\delta \colon \che{ M} \rightarrow \che{ N}$ such that for any $w \in N$ and any $\xi \in \che{ M}$ it implies
\[
\pairx{ w}{ \delta( \xi)}= \pairx{ p( \gamma( w))}{ \xi}.
\]
By the commutative diagram
\[
\squarerd{ M}{ \alpha}{ \iota_{ M}}{ N}{ \gamma}{ \che{( \che{ M})}}{ i}{ \cheloc{ \che{ M}}}
\]
for any $v \in M$ we have
\[
\pairx{ \alpha( v)}{ \delta( \xi)}= \pairx{ v}{ \xi}.
\]
So $\che{ \alpha} \circ \delta= \xxid_{ \che{ M}}$. So we have
\[
\xxdim( \che{ M}) \leq \xxdim( \che{ N})= n.
\]
By Lemma \ref{Prop:FA9} and Proposition \ref{Prop:FA8}, there is an injection from the set of the extremal rays of $M$ to the set of the extremal rays of $\che{ M}$. So we have $\xxdim( M) \leq n$.
\end {proof}

\begin {Ex}
There is an example of straight submodule $M \subset \rrt^2$ that is not finitely generated. Let
\[
M= \{ ( a, b) \in \rrt^2 \setmid b \not= \mx \} \cup \{ \mx \}.
\]
$M$ is a submodule of $\rrt^2$. $M$ is straight, because it is lattice-preserving.
\end {Ex}

\begin {Ex} \label {Ex:FA2}
There is an example of extremally generated submodule $M \subset \rrt^3$ that is not finitely generated. Let
\[
M= \left\{
\begin {array}[c]{rcl}
( a, b, c) \in \rrt^3 & \setmid & ( -1) \dx a \px c \lqx b, \\
&& 2 b \lqx a \dx c
\end {array}
\right\}.
\]
$M$ is a submodule of $\rrt^3$ (Example \ref{Ex:FA3}). For $0 \lqx t \lqx 1$, let
\[
e( t)=( 2t, t, 0) \in M.
\]
$e( t)$ is extremal. (Proposition \ref{Prop:FA4}. Indeed $e( t)$ is a minimal element of the subset
\[
S_t= \{( a, b, c) \in M \setmid b= t \}.
\]
So it is extremal.) So $M$ is not finitely generated. $\{ e( t) \setmid 0 \lqx t \lqx 1 \}$ is a basis of $M$. Indeed, for any $( a, b, c) \in M$,
\[
( a, b, c)= c \dx e( b \sx c) \px( 2 b \sx a) \dx e( a \sx b).
\]
$M$ is not straight. Indeed, let
\begin {align*}
v_1 & =( \frac{ 1}{ 2}, \frac{ 1}{ 2}, \frac{ 1}{ 2}), \\
v_2 & =( 1, \frac{ 1}{ 2}, 0), \\
w & =( 1, 0, 0).
\end {align*}
Then we have
\begin {align*}
\infx{ M}{ v_1}{ v_2} & =( \frac{ 1}{ 2}, \frac{ 1}{ 4}, 0), \\
\infx{ M}{ v_1}{ v_2} \px w & =( 1, \frac{ 1}{ 4}, 0), \\
\infx{ M}{ v_1 \px w}{ v_2 \px w} & =( 1, \frac{ 1}{ 2}, 0).
\end {align*}
So $M$ is not straight.
\end {Ex}

\subsection {Free modules}

Let $k$ be a totally ordered tropical semifield. Let $F= k^n$ be the free module with the basis $\{ e_1, \ldots, e_n \}$. Let $F^*$ be the set of the linear combinations of $\{ e_1, \ldots, e_n \}$ with coefficients in $k ^*= k \setminus \{ \mx \}$. Let $\{ \che{ e}_1, \ldots, \che{ e}_n \}$ be the dual basis in $\che{ F}$. We have a bijective map
\[
\psi \colon F^* \rmap ( \che{ F})^*
\]
defined by
\[
\psi( a_1 \dx e_1 \px \cdots \px a_n \dx e_n)=( \sx a_1) \dx \che{ e}_1 \px \cdots \px( \sx a_n) \dx \che{ e}_n.
\]
For $v, w \in F^*$, the condition $v \lqx w$ is fulfilled if and only if
\[
\pairx{ v}{ \psi( w)} \lqx 0.
\]
For $w \in F^*$ and $1 \leq i \leq n$, let
\[
M( w, i)= \{ v \in F \setmid \setforall j, \pairx{ v}{ \che{ e}_j} \dx \pairx{ e_j}{ \psi( w)} \gqx \pairx{ v}{ \che{ e}_i} \dx \pairx{ e_i}{ \psi( w)} \}.
\]
$M( w, i)$ is a submodule of $F$ (Example \ref{Ex:FA3}). It is easy to see that $M( w, i)$ is lattice-preserving in $F$, i.e. the inclusion $M( w, i) \rightarrow F$ preserves the infimum of any two elements. For $\eta \in \che{ F}$ and $1 \leq i \leq n$, let
\begin {align*}
N( \eta, i) & = \{ v \in F \setmid \pairx{ v}{ \eta}= \pairx{ v}{ \che{ e}_i} \dx \pairx{ e_i}{ \eta} \} \\
& = \{ v \in F \setmid \setforall j, \pairx{ v}{ \che{ e}_j} \dx \pairx{ e_j}{ \eta} \lqx \pairx{ v}{ \che{ e}_i} \dx \pairx{ e_i}{ \eta} \}.
\end {align*}
$N( \eta, i)$ is also a lattice-preserving submodule of $F$.

\begin {Prop} \label {Prop:FF1}
Let $M$ be a submodule of $F$ with a basis $\{ w_1, \ldots, w_r \}$. Suppose that $w_h \in F^*$ ($1 \leq h \leq r$). Then the following are equivalent.
\begin {enumerate}
\renumi{\conditionenumi}
\item $M$ is lattice-preserving in $F$.
\item For any $i \in \{ 1, \ldots, n \}$, there is the minimum element of $M \cap V_i$, where
\[
V_i= \{ v \in F \setmid \pairx{ v}{ \che{ e}_i}= 0 \}.
\]
\item There is a surjective map
\[
s \colon \{ 1, \ldots, n \} \rmap \{ 1, \ldots, r \}
\]
such that
\[
M= \bigcap_{ 1 \leq i \leq n} M( w_{ s( i)}, i).
\]
\item There is a surjective map
\[
s \colon \{ 1, \ldots, n \} \rmap \{ 1, \ldots, r \}
\]
such that
\[
M= \bigcap_{ 1 \leq i \leq n} N( \eta_{ s( i)}, i),
\]
where $\eta_h$ is the dual element of $w_h$.
\end {enumerate}
\end {Prop}

\begin {proof}
\condition{iii} $\then$ \condition{i} and \condition{iv} $\then$ \condition{i} are easy.
\par
Since $\che{ F}$ is also a free module, for $\eta \in ( \che{ F})^*$ and $1 \leq i \leq n$ we have the lattice-preserving submodule $M( \eta, i)$ of $\che{ F}$. The bijective map
\[
\psi \colon F^* \rmap ( \che{ F})^*
\]
induces bijective maps
\[
\psi' \colon M \setminus \{ \mx \} \rmap \che{ M} \setminus \{ \mx \},
\]
\[
\psi'' \colon N( \eta, i) \setminus \{ \mx \} \rmap M( \eta, i) \setminus \{ \mx \}.
\]
So we have only to prove that conditions \condition{i}, \condition{ii}, \condition{iii} are equivalent.
\par
\condition{i} $\then$ \condition{ii}. Let
\[
v_i= \inf_{ F} \{ \sx a_{ h, i} \dx w_h \setmid 1 \leq h \leq r \},
\]
where
\[
a_{ h, i}= \pairx{ w_h}{ \che{ e}_i}.
\]
Then $v_i$ is the minimum element of $M \cap V_i$.
\par
\condition{ii} $\then$ \condition{iii}. Let $v_i$ be the minimum element of $M \cap V_i$. $v_i$ is an extremal element of $M$. The extremal ray $k \dx v_i$ is generated by an element of the basis $\{ w_1, \ldots, w_r \}$ (Proposition \ref{Prop:LK6}). There is a number $s( i)$ such that $k \dx v_i= k \dx w_{ s( i)}$.
\par
We show that $s$ is surjective. For $h \in \{ 1, \ldots, r \}$, we have
\[
w_h= a_{h, 1} \dx v_1 \px \cdots \px a_{ h, n} \dx v_n.
\]
Since $w_h$ is extremal (Proposition \ref{Prop:LK6}), there is a number $i$ such that
\[
w_h= a_{ h, i} \dx v_i.
\]
So we have $h= s( i)$.
\par
We show the equality
\[
M= \bigcap_{ 1 \leq i \leq n} M( w_{ s( i)}, i).
\]
For $v \in F$, let $x_i= \pairx{ v}{ \che{ e}_i}$. The condition $v \in M( w_{ s( i)}, i)$ is fulfilled if and only if for any $j$ it implies
\[
\sx a_{ s( i), j} \dx x_j \gqx \sx a_{ s( i), i} \dx x_i.
\]
For $1 \leq h \leq r$ and $1 \leq i \leq n$, we have
\begin {align*}
\sx a_{ s( i), i} \dx w_{ s( i)} & = v_i \\
& \lqx \sx a_{ h, i} \dx w_{ h}.
\end {align*}
For $1 \leq j \leq n$, we have
\[
\sx a_{ s( i), i} \dx a_{ s( i), j} \lqx \sx a_{ h, i} \dx a_{ h, j}.
\]
It means $w_h \in M( w_{ s( i)}, i)$. Since $M$ is generated by $\{ w_1, \ldots, w_r \}$, we have
\[
M \subset \bigcap_{ 1 \leq i \leq n} M( w_{ s( i)}, i).
\]
Let $v$ be an element of right side. Then
\[
v= \bigpx_{ i} \sx a_{ s( i), i} \dx x_i \dx w_{ s( i)}.
\]
(Indeed,
\begin {align*}
\pairx{ v}{ \che{ e}_i} & = x_i \\
& = \pairx{ \sx a_{ s( i), i} \dx x_i \dx w_{ s( i)}}{ \che{ e}_i}.
\end {align*}
So
\[
v \lqx \bigpx_{ i} \sx a_{ s( i), i} \dx x_i \dx w_{ s( i)}.
\]
The converse is easy.) So we have $v \in M$.
\end {proof}

\section {Polytopes in a tropical projective space} \label {Sec:FB3}

Let $F= \rrt^{ n+1}$ be the free module with coordinates $( x_1, \ldots, x_{ n+1})$ over $\rrt= \rr \cup \{ \mx \}$. Let $F^*= \rr^{ n+1}$.

\begin {Prop} \label {Prop:FF2}
Let $M$ be a submodule of $F$ generated by finitely many elements of $F^*$. Then the following are equivalent.
\begin {enumerate}
\renumi{\conditionenumi}
\item $M$ is lattice-preserving in $F$.
\item $M \setminus \{ \mx \}$ is a real convex subset of $\rr^{ n+1}$.
\end {enumerate}
\end {Prop}

\begin {proof}
\condition{i} $\then$ \condition{ii}. By Proposition \ref{Prop:FF1}, $M$ is defined by inequalities
\[
x_j \gqx x_i- c_{ i, j} \quad ( i, j \in \{ 1, \ldots, n+1 \})
\]
for some $c_{ i, j} \in \rr$. So $M \setminus \{ \mx \}$ is real convex.
\par
\condition{ii} $\then$ \condition{i}. Let $\pi_1 \colon F \rightarrow \rrt^n$ and $\pi_2 \colon F \rightarrow \rrt$ be projections defined as follows.
\begin {align*}
\pi_1( x_1, \ldots, x_{ n+1}) & =( x_1, \ldots, x_n), \\
\pi_2( x_1, \ldots, x_{ n+1}) & = x_{ n+1}.
\end {align*}
For $a \in \rr$, let $N_i( a) \subset F$ be the submodule defined as follows.
\[
N_i( a)= \{ v=( x_1, \ldots, x_{ n+1}) \in F \setmid x_{ n+1}= x_i+ a \}.
\]
By induction on $n$, we may assume that modules $\pi_1( M)$, $\pi_2( M)$, $M \cap N_i( a)$ are lattice-preserving. Suppose that $M$ is not lattice-preserving. By Proposition \ref{Prop:FF1}, there is a number $i$ such that there is no minimum element of $M \cap V_i$, where
\[
V_i= \{ v=( x_1, \ldots, x_{ n+1}) \in F \setmid x_i= 0 \}.
\]
We may assume $i \leq n$. Let $w_1, w_2$ be minimal elements of $M \cap V_i$ such that $\pi_1( w_1)$ is the minimum element of $\pi_1( M \cap V_i)$ and that $\pi_2( w_2)$ is the minimum element of $\pi_2( M \cap V_i)$. Let $a \in \rr$ be an element such that
\[
\pi_2( w_2) \ltx a \ltx \pi_2( w_1).
\]
There is the minimum element $v( a)$ of $M \cap N_i( a) \cap V_i$. Since $M \cap V_i$ is real convex, $v( a)$ is a minimal element of $M \cap V_i$. (Let $v' \in M \cap V_i$ be an element such that $v' \ltx v( a)$. The real line segment combining $v'$ and $w_1$ contains an element $v'' \in M \cap N_i( a) \cap V_i$ such that $v'' \not= v'$. Since $\pi_1( w_1) \ltx \pi_1( v') \lqx \pi_1( v( a))$, we have $v'' \ltx v( a)$. ) So $M$ has infinitely many extremal rays, which is contradiction.
\end {proof}

Let
\[
\ph \colon \rrt^{ n+1} \setminus \{ \mx \} \rmap \rrt \pp^n
\]
be the canonical projection to the tropical projective space $\rrt \pp^n$. We identify $\ph( \rr^{ n+1})$ with $\rr^n$. A subset $P \subset \rrt \pp^n$ is said to be tropically convex if the subset
\[
M= \inv{ \ph}( P) \cup \{ \mx \} \subset \rrt^{ n+1}
\]
is a submodule. A subset $P \subset \rrt \pp^n$ is said to be a tropical polytope if it is the tropically convex hull of finitely many points of $\rr^n$.

\begin {proof}[Proof of Theorem \ref{Prop:FE5}]
\statement{1} Suppose that $P$ is a polytrope. Then $P$ is real convex. By Proposition \ref{Prop:FF2}, $M$ is lattice-preserving in $\rrt^{ n+1}$. So $M$ is straight. By Theorem \ref{Prop:FE3}, $M$ is reflexive.
\par
\statement{2} Suppose that $M$ is straight reflexive. Let $\{ v_1, \ldots, v_r \}$ be a basis of $M$. By Theorem \ref{Prop:FE2}, we have $r \leq n+1$. Let $p_i= \ph( v_i)$. Then $P$ is the tropically convex hull of $\{ p_1, \ldots, p_r \}$.
\end {proof}

\section {Square matrices over a tropical semifield} \label {Sec:FB4}

Let $k$ be a totally ordered rational tropical semifield. A square matrix of order $n$ over $k$ is a homomorphism $A \colon k^n \rightarrow k^n$. Let $\{ e_1, \ldots, e_n \}$ be the basis of $k^n$. The coefficient $\pairx{ A \dx e_j}{ \che{ e_i}}$ is simply written as $A_{ i j}$. Let $E_n \colon k^n \rightarrow k^n$ be the identity.
\par
Let $\Delta( A)$, $\overline{ \Delta}( A)$ be square matrices of order $n$ defined as follows.
\[
\Delta( A)_{ i j}= \delta_{ i j} \dx A_{ i j},
\]
\[
\overline{ \Delta}( A)_{ i j}= \overline{ \delta}_{ i j} \dx A_{ i j},
\]
where
\[
\delta_{ i j}= \begin {cases}
0 & \text{if $i= j$} \\
\mx & \text{if $i \not= j$}
\end {cases}
\]
\[
\overline{ \delta}_{ i j}= \begin {cases}
\mx & \text{if $i= j$} \\
0 & \text{if $i \not= j$}.
\end {cases}
\]
The determinant $\xxdet( A)$ is the sum of elements $A_{ 1 s( 1)} \dx \cdots \dx A_{ n s( n)}$ for all permutations $s \in S( n)$.

\begin {Lem} \label {Prop:FF3}
Let $A$ be a square matrix of order $n$ over $k$. Suppose that $\Delta( A)= E_n$ and $\xxdet( A)= 0$. Then $A^{ \dx n}= A^{ \dx n-1}$.
\end {Lem}

\begin {proof}
Since $E_n \lqx A$, we have $A^{ \dx r} \lqx A^{ \dx r+1}$ for any $r \geq 0$. $( A^{ \dx n})_{ i j}$ is the sum of elements
\[
b= A_{ h( 0) h( 1)} \dx A_{ h( 1) h( 2)} \dx \cdots \dx A_{ h( n-1) h( n)}
\]
for all maps $h \colon \{ 0, \ldots, n \} \rightarrow \{ 1, \ldots, n \}$ such that $h( 0)= i$ and $h( n)= j$. $h$ is not injective. So there are numbers $l, m$ and a cyclic permutation $s \in S( n)$ such that
\[
s \colon h( l) \mapsto h( l+1) \mapsto \cdots \mapsto h( m-1) \mapsto h( m)= h( l).
\]
Since $\Delta( A)= E_n$, we have
\[
A_{ h( l) h( l+1)} \dx \cdots \dx A_{ h( m-1) h( m)} \lqx \xxdet( A).
\]
So we have
\[
A^{ \dx n} \lqx \xxdet( A) \dx A^{ \dx n-1}.
\]
Since $\xxdet( A)= 0$, we have the conclusion.
\end {proof}

\begin {Lem} \label {Prop:FF4}
Let $A$ be a square matrix of order $n$ over $k$. Then either \condition{i} or \condition{ii} is fulfilled.
\begin {enumerate}
\renumi{\conditionenumi}
\item There are an element $v \in( k \setminus \{ \mx \})^n$ and an element $\e \gtx 0$ such that
\[
( A \px \e \dx \overline{ \Delta}( A)) \dx v= \Delta( A) \dx v.
\]
\item There is an element $v \in k^n \setminus \{ \mx \}$ such that
\[
A \dx v= \overline{ \Delta}( A) \dx v.
\]
\end {enumerate}
\end {Lem}

\begin {proof}
Let $e( A)$ be the sum of elements $A_{ 1 s( 1)} \dx \cdots \dx A_{ n s( n)}$ for all $s \in S( n) \setminus \{ \xxid \}$. Let
\[
c( A)= \xxdet( \Delta( A))= A_{ 1 1} \dx \cdots \dx A_{ n n}.
\]
\par
We show that the condition \condition{i} is fulfilled if $e( A) \ltx c( A)$. Replacing $A$ by $ \sx( \Delta( A)) \dx A$, we may assume $\Delta( A)= E_n$. There is an element $\e \in k$ such that $\e \gtx 0$ and
\[
e( A) \dx n \e \lqx c( A).
\]
Let
\[
B= A \px \e \dx \overline{ \Delta}( A).
\]
Then we have $e( B) \lqx c( B)$. By Lemma \ref{Prop:FF3}, we have $B^{ \dx n}= B^{ \dx n-1}$. Let $w \in( k \setminus \{ \mx \})^n$ be any element. Let $v= B^{ \dx n-1} \dx w$. Then we have $B \dx v= v$.
\par
We show that the condition \condition{ii} is fulfilled if $c( A) \lqx e( A)$. We may assume $\Delta( A)= E_n$. (If $A_{ i i}= \mx$, then the element $v= e_i$ satisfies the conclusion.) There is a cyclic permutation $s \in S( n) \setminus \{ \xxid \}$ and a map $h \colon \{ 0, \ldots, l \} \rightarrow \{ 1, \ldots, n \}$ such that
\[
s \colon h( 0) \mapsto h( 1) \mapsto \cdots \mapsto h( l-1) \mapsto h( l)= h( 0),
\]
\[
A_{ h( 0) h( 1)} \dx \cdots \dx A_{ h( l-1) h( l)} \gqx 0.
\]
Let
\[
v= \bigpx_{ 1 \leq m \leq l}( A_{ h( m) h( m+1)} \dx \cdots  \dx A_{ h( l-1) h( l)}) \dx e_{ h( m)}.
\]
Then
\[
\overline{ \Delta}( A) \dx v \gqx v.
\]
So we have the conclusion.
\end {proof}

\section {Tropical curves} \label {Sec:FB5}

Let $A= \rrt[x_1, - x_1, \ldots, x_n, - x_n]$ be the semiring of Laurent polynomials over $\rrt= \rr \cup \{ \mx \}$ (where $-x_i$ means $\sx x_i$). Let
\[
f= \bigpx_{ i_1 \ldots i_n \in \zz} c_{ i_1 \ldots i_n} \dx i_1 x_1 \dx \cdots \dx i_n x_n
\]
be any element of $A$. The induced map
\[
\begin {array}[c]{rcl}
f \colon \rr^n & \rmap & \rrt \\
( a_1, \ldots, a_n) & \mapsto & f( a_1, \ldots, a_n)
\end {array}
\]
is said to be a Laurent polynomial function over $\rrt$. If $f$ is a monomial, then $f$ is a $\zz$-affine function, i.e. there are $c \in \rr$ and $i_1, \ldots, i_n \in \zz$ such that
\[
f= c+ i_1 x_1+ \cdots + i_n x_n.
\]
In general case, $f$ is the supremum of finitely many $\zz$-affine functions, which is a locally convex piecewise-$\zz$-affine function.
\par
Let $\Gamma_n \subset \rr^n$ be the subset defined as follows.
\[
\Gamma_n= E_0 \cup E_1 \cup \cdots \cup E_n,
\]
\[
E_0= \{( a_1, \ldots, a_n) \in \rr^n \setmid \setforall i, \setforall j, a_i= a_j \gqx 0 \},
\]
for $1 \leq i \leq n$,
\[
E_i= \{( a_1, \ldots, a_n) \in \rr^n \setmid a_i \lqx 0, \setforall j \not= i, a_j= 0 \}.
\]
$\Gamma_n$ has a $(n+1)$-valent vertex $P=( 0, \ldots, 0)$. Also $\Gamma_n$ is equipped with Euclidean topology on $\rr^n$.

\begin {Def}
A function $f \colon \Gamma_n \rightarrow \rrt$ is {\em regular} if it is induced by a locally Laurent polynomial function $f \colon \rr^n \rightarrow \rrt$.
\end {Def}

Let $\ooo_{ \Gamma_n}$ be the sheaf of the regular functions on $\Gamma_n$. $\ooo_{ \Gamma_n}$ is a sheaf of semirings. Let $R$ be the stalk of $\ooo_{ \Gamma_n}$ at the vertex $P$.

\begin {Prop} \label {Prop:FF5}
Let $f \in R \setminus \{ \mx \}$ be any element. Then there are a unique number $r \in \zz_{ \geq 0}$ and a unique Laurent monomial $h \in R$ such that
\[
f= h \dx r ( x_1 \px 0).
\]
\end {Prop}

\begin {proof}
$f$ is the sum of Laurent monomials $f_1, \ldots, f_m$. If $f_j( P) \ltx f( P)$, then $f_j \ltx f$ on a neighborhood of $P$. So we may assume $f_j( P)= f( P)$. Then $f$ is $\zz$-affine on $E_i$ ($1 \leq i \leq n$). So there is $a_i \in \zz$ such that $f= f( P) \dx a_i x_i$ on $E_i$. Let
\[
h= f( P) \dx a_1 x_1 \dx \cdots \dx a_n x_n.
\]
Then $f= h$ on $E_1 \cup \cdots \cup E_n$. $f \sx h$ is the sum of monomials $g_1, \ldots, g_m$ such that $g_j( P)= 0$. There are $b_{ i j} \in \zz_{ \geq 0}$ such that
\[
g_j= b_{ 1 j} x_1 \dx \cdots \dx b_{ n j} x_n.
\]
Then
\[
g_j=( b_{ 1 j}+ \cdots + b_{ n j}) x_1
\]
on $E_0$. So we have $f \sx h= r x_1$ on $E_0$, where
\[
r= \bigpx_{ 1 \leq j \leq m} \sum_{ 1 \leq i \leq n} b_{ i j}.
\]
\end {proof}

The number $r$ in the above statement is called the order of $f$ at $P$, and denoted by $\xxord( f, P)$.
\par
For $0 \leq i \leq n$, let $X_i f$ be the partial differential of $f$ at $P$ with direction $E_i$. (i.e. $X_i f= a$ if and only if $f= f( P)- a x_i$ on $E_i$ ($1 \leq i \leq n$). $X_0 f= a$ if and only if $f= f( P)+ a x_1$ on $E_0$.)

\begin {Prop} \label {Prop:FF7}
Let $f \in R \setminus \{ \mx \}$ be any element. Then
\[
\xxord( f, P)= \sum_{ 0 \leq i \leq n} X_i f.
\]
\end {Prop}

\begin {proof}
Let $h$ be a Laurent monomial written as follows.
\[
h= c \dx a_1 x_1 \dx \cdots \dx a_n x_n.
\]
Then
\[
X_i h= -a_i \quad ( 1 \leq i \leq n),
\]
\[
X_0 h= a_1+ \cdots + a_n.
\]
So we have
\[
\sum_{ 0 \leq i \leq n} X_i h= 0.
\]
Also we have
\[
\sum_{ 0 \leq i \leq n} X_i( x_1 \px 0)= 1.
\]
So we have the conclusion.
\end {proof}

\begin {Prop} \label {Prop:FF6}
Let $f, g \in R \setminus \{ \mx \}$ be any elements.
\begin {enumerate}
\renumi{\statementenumi}
\item $\xxord( f \dx g, P)= \xxord( f, P)+ \xxord( g, P)$.
\item If $g( P) \lqx f( P)$, then $\xxord( f, P) \lqx \xxord( f \px g, P)$.
\end {enumerate}
\end {Prop}

\begin {proof}
\statement{1} is easy.
\par
\statement{2} If $g( P) \ltx f( P)$, then $f \px g= f$. So we may assume $g( P)= f( P)$. Then we have
\[
X_i( f \px g)= X_i f \px X_i g.
\]
By Proposition \ref{Prop:FF7}, we have the conclusion.
\end {proof}

A function $f \colon \Gamma_n \rightarrow \rrt$ is said to be rational if locally
\[
f= g_1- g_2= g_1 \sx g_2
\]
for regular functions $g_1, g_2$. By Proposition \ref{Prop:FF5}, there is a number $m \geq 0$ such that the function $m( x_1 \px 0) \dx f$ is regular at $P$. The order of $f$ at $P$ is defined as follows.
\[
\xxord( f, P)= \xxord( m( x_1 \px 0) \dx f, P)- m.
\]
\par
Let $Q \in \Gamma_n$ be a point such that $Q \not= P$. Then a neighborhood of $Q$ is embedded in $\Gamma_1= \rr$. So we can define the order of $f$ at $Q$ similarly.

\begin {Def}
$(C, \ooo_C)$ is a {\em tropical curve} if for any $P \in C$ there are a neighborhood $U$ of $P$ and a number $n \geq 1$ such that $(U, \ooo_U)$ is embedded in $( \Gamma_n, \ooo_{ \Gamma_n})$.
\end {Def}

A divisor $D$ on a tropical curve $C$ is an element of the free abelian group $\xxDiv( C)$ generated by all the points of $C$. For a rational function $f \colon C \rightarrow \rrt$, the divisor $( f) \in \xxDiv( C)$ is defined as follows.
\[
( f)= \sum_{ P \in C} \xxord( f, P) P.
\]
$f$ is said to be a section of $D$ if either $f= \mx$ or $( f)+ D \geq 0$. Let $\ooo_C( D)$ be the sheaf of the sections of $D$.

\begin {Prop} \label {Prop:FF8}
The set $M= H^0( C, \ooo_C( D))$ is a $\rrt$-module.
\end {Prop}

\begin {proof}
Let $f, g \in M \setminus \{ \mx \}$ be any elements. By Proposition \ref{Prop:FF6}, for $P \in C$ we have
\[
\xxord( f \px g, P) \gqx \min \{ \xxord( f, P), \xxord( g, P) \}.
\]
So
\[
( f \px g)+ D \gqx \xxinf_{ \xxDiv( C)} \{( f),( g) \}+ D \gqx 0.
\]
So we have $f \px g \in M$.
\end {proof}

Recall that
\[
r( D)= \max \{ r \in \zz_{ \geq -1} \setmid U( D, r) = \emptyset \}.
\]

\begin {proof}[Proof of Theorem \ref{Prop:FE7}]
Note that $r( D)= s( D)-1$, where
\[
s( D)= \min \{ r \in \zz_{ \geq 0} \setmid U( D, r) \not= \emptyset \}.
\]
Let $m= s( D)$. We show that there is a straight reflexive submodule $N \subset M= H^0( C, \ooo_{ C}( D))$ with dimension $m$. Let $P_1, \ldots, P_m \in C$ be points such that
\[
H^0( C, \ooo_{ C}( D- E))= \mx,
\]
where
\[
E= P_1+ \cdots + P_m.
\]
There is an element
\[
f_i \in H^0( C, \ooo_{ C}( D- E+ P_i))
\]
such that $f_i \not= \mx$. Let
\[
\alpha \colon \rrt^m \rmap M
\]
be the homomorphism defined by $\alpha( e_i)= f_i$. Let
\[
\beta \colon M \rmap \rrt^m
\]
be the homomorphism defined by
\[
\beta( g)= g( P_1) \dx e_1 \px \cdots \px g( P_m) \dx e_m.
\]
Let $A$ be the square matrix induced by $\beta \circ \alpha \colon \rrt^m \rightarrow \rrt^m$.
\par
Now we suppose that there is an element
\[
v= a_1 \dx e_1 \px \cdots \px a_m \dx e_m \in \rrt^m \setminus \{ \mx \}
\]
such that $A \dx v= \overline{ \Delta}( A) \dx v$. Then there is a map $h \colon \{ 1, \ldots, m \} \rightarrow \{ 1, \ldots, m \}$ such that $h( i) \not= i$ and
\[
\alpha( v)( P_i)= a_{ h( i)} \dx f_{ h( i)}( P_i).
\]
Then
\[
\xxord( \alpha( v), P_i) \gqx \xxord( f_{ h( i)}, P_i)
\]
(Proposition \ref{Prop:FF6}). So $\alpha( v)$ is a section of $D- E$ such that $\alpha( v) \not= \mx$, which is contradiction.
\par
So there is no element $v \in \rrt^m \setminus \{ \mx \}$ such that $A \dx v= \overline{ \Delta}( A) \dx v$. By Lemma \ref{Prop:FF4}, there are an element $v \in \rr^m$ and an element $\e \gtx 0$ such that
\[
( A \px \e \dx \overline{ \Delta}( A)) \dx v= \Delta( A) \dx v.
\]
Let $L( v, \e) \subset \rrt^m$ be the submodule defined as follows.
\[
L( v, \e)= \rrt \dx \{ w \in \rrt^m \setmid v \lqx w \lqx \e \dx v \}.
\]
$L( v, \e)$ is a straight reflexive $\rrt$-module with dimension $m$. We have
\[
A \res{ L( v, \e)}= \Delta( A) \res{ L( v, \e)}.
\]
So $\alpha$ is injective on $L( v, \e)$. The image $N= \alpha( L( v, \e))$ is a submodule of $M$ such that $N \iso L( v, \e)$.
\end {proof}

\begin {Ex} \label {Ex:FB1}
The mapping $D \mapsto r( D)$ is not an invariant of a $\rrt$-module. We show that there are tropical curves $C, C'$ and divisors $D, D'$ such that
\[
H^0( C, \ooo_C( D)) \iso H^0( C', \ooo_{C'}( D')),
\]
\[
r( D) \not= r( D').
\]
\par
Let $C$ be a tropical curve with genus $1$ with a vertex $V$ and an edge $E$. Let $P$ be an interior point of $E$. Let $D= V+ P$. Then $H^0( C, \ooo_C( D))$ is isomorphic to the submodule of $\rrt^2$ generated by $( 0, 0)$ and $( 0, \frac{ a}{ 2})$, where $a$ is the lattice length of $E$. We have $r( D)= 1$.
\par
Let $C'$ be a tropical curve with genus $2$ with vertices $V_1, V_2$ and edges $E_1, E_2, E_3$ such that the boundary of $E_i$ is $\{ V_1, V_2 \}$ ($1 \leq i \leq 3$). Let $P$ be an interior point of $E_1$. Let $D'= V_1+ P$. Then for any interior point $Q$ of $E_2 \cup E_3$ we have
\[
H^0( C', \ooo_{C'}( D'- Q))= \mx.
\]
So $H^0( C', \ooo_{C'}( D'))$ is isomorphic to the submodule of $\rrt^2$ generated by $( 0, 0)$ and $( 0, \frac{ b}{ 2})$, where $b$ is the lattice length of the path from $V_1$ to $P$ contained in $E_1$. We have $r( D')= 0$. In the case of $a= b$, the required condition is fulfilled.
\end {Ex}

\section {Tropical plane curves} \label {Sec:B1}

\subsection {Tropicalization}

It is well known that some example of tropical curve is given by tropicalization of a family of affine complex curves.
\par
First, we define tropical plane curves. Let $f \in \rrt[x, - x, y, - y]$ be a Laurent polynomial over $\rrt= \rr \cup \{ \mx \}$. The subset
\[
V( f)= \{ (a, b) \in \rr^2 \setmid -f \text{ is not locally convex at } (a, b) \}
\]
is called the algebraic subset defined by $f$. The morphism $C_f \rightarrow \rr^2$ parametrizing $V( f)$ with a tropical curve $C_f$ is called the tropical plane curve defined by $f$. The genus of $C_f$ is defined to be the first Betti number $b_1( C_f)$.
\par
A tropical plane curve is a dequantization of complex amoebas in following way. For $t > 1$, let
\[
\aaa_t \colon ( \cc^{ \times})^2 \rmap \rr^2
\]
be the homomorphism of groups defined by
\[
\aaa_t( a, b)=( \frac{ \log \abso{ a}}{ \log( t)}, \frac{ \log \abso{ b}}{ \log( t)}).
\]
$\aaa_t$ is called the complex amoeba map. Let
\[
g_t \in \cc[z_1, \inv{ z_1}, z_2, \inv{ z_2}] \quad (t > 1)
\]
be a family of complex Laurent polynomials such that each coefficient is a Laurent polynomial of $\inv{ t}$. This family is written as an element of a valuation field $K$. We use the group algebra $K= \cc[[ \rr]]$ of power series defined by the group $\rr$. The indeterminate is denoted by $\inv{ t}$, and the valuation is defined to be the maximum index of $t$ multiplied by $-1$. So, $\xxval( t^a)= -a$. The family $\{ g_t \setmid t > 1 \}$ is written as an element
\[
g \in K[z_1, \inv{ z_1}, z_2, \inv{ z_2}].
\]
The amoeba map over $K$
\[
\aaa \colon ( K^{ \times})^2 \rmap \rr^2
\]
is defined as follows.
\[
\aaa( a, b)=( - \xxval( a), - \xxval( b)).
\]
The affine curve $V( g) \subset ( K^{ \times})^2$ is the family of affine complex curves $V( g_t) \subset ( \cc^{ \times})^2$. Taking $t \rightarrow + \infty$, the family of complex amoebas $\aaa_t( V( g_t))$ converges to the amoeba $\aaa( V( g))$ over $K$. Also, the amoeba over $K$ is the algebraic subset defined by a tropical Laurent polynomial. Let
\[
\aaa \colon K[z_1, \inv{ z_1}, z_2, \inv{ z_2}] \rmap \rrt[x, - x, y, - y]
\]
be the map defined as follows.
\[
\aaa( g)= f,
\]
\[
g= \sum_{ i, j \in \zz} c_{ i j} z_1^i z_2^j,
\]
\[
f= \bigpx_{ i, j \in \zz} - \xxval( c_{ i j}) \dx i x \dx j y.
\]
Then we have
\[
\aaa( V( g))= V( f).
\]
This construction is called the tropicalization of a family of affine complex curves.

\subsection {Examples}

\begin {Ex}
For $a, b, c \in \cc^{ \times}$, let
\[
g= a+ b z_1+ c z_2.
\]
Then
\[
f= \aaa( g)= 0 \px x \px y.
\]
The tropical plane curve $C_f$ is said to be a tropical projective line. We have $b_1( C_f)= 0$.
\end {Ex}

\begin {Ex}
For $r, s \in \nn$ and $a_i, b_j \in \rr$, let
\[
f= f_1 \dx f_2,
\]
\[
f_1= a_0 \px a_1 \dx x \px a_2 \dx 2 x \px \cdots \px a_r \dx r x,
\]
\[
f_2= b_0 \px b_1 \dx y \px b_2 \dx 2 y \px \cdots \px b_s \dx s y.
\]
Assume that
\[
2 a_i \gtx a_{ i-1}+ a_{ i+1},
\]
\[
2 b_j \gtx b_{ j-1}+ b_{ j+1}.
\]
Then $b_1( C_f)=( r-1)( s-1)$.
\end {Ex}

\end {document}